# UNIFORM LIMIT THEOREMS FOR WAVELET DENSITY ESTIMATORS

By Evarist Giné and Richard Nickl

*University of Connecticut and University of Cambridge*

Let $p_n(y) = \sum_k \hat{\alpha}_k \phi(y-k) + \sum_{l=0}^{j_n-1} \sum_k \hat{\beta}_{lk} 2^{l/2} \psi(2^l y - k)$ be the linear wavelet density estimator, where $\phi, \psi$ are a father and a mother wavelet (with compact support), $\hat{\alpha}_k, \hat{\beta}_{lk}$ are the empirical wavelet coefficients based on an i.i.d. sample of random variables distributed according to a density $p_0$ on $\mathbb{R}$, and $j_n \in \mathbb{Z}$, $j_n \nearrow \infty$. Several uniform limit theorems are proved: First, the almost sure rate of convergence of $\sup_{y \in \mathbb{R}} |p_n(y) - Ep_n(y)|$ is obtained, and a law of the logarithm for a suitably scaled version of this quantity is established. This implies that $\sup_{y \in \mathbb{R}} |p_n(y) - p_0(y)|$ attains the optimal almost sure rate of convergence for estimating $p_0$, if $j_n$ is suitably chosen. Second, a uniform central limit theorem as well as strong invariance principles for the distribution function of $p_n$, that is, for the stochastic processes $\sqrt{n}(F_n^W(s) - F(s)) = \sqrt{n} \int_{-\infty}^s (p_n - p_0), s \in \mathbb{R}$, are proved; and more generally, uniform central limit theorems for the processes $\sqrt{n} \int (p_n - p_0) f, f \in \mathcal{F}$, for other Donsker classes $\mathcal{F}$ of interest are considered. As a statistical application, it is shown that essentially the same limit theorems can be obtained for the hard thresholding wavelet estimator introduced by Donoho et al. [*Ann. Statist.* **24** (1996) 508–539].

**1. Introduction.** Let $X, X_1, \ldots, X_n$ be independent identically distributed real-valued random variables with absolutely continuous law $P$ and density $p_0$, and denote by $P_n$ the usual empirical measure induced by the sample. If $\phi$ is a bounded and compactly supported father wavelet (scaling function) and $\psi$ an associated mother wavelet, the (linear) wavelet density estimator of $p_0$ is given by

$$p_n(y) = \sum_{k \in \mathbb{Z}} \hat{\alpha}_{j_n k} 2^{j_n/2} \phi(2^{j_n} y - k)$$

(1)
$$= \sum_{k \in \mathbb{Z}} \hat{\alpha}_{0k} \phi(y-k) + \sum_{l=0}^{j_n-1} \sum_{k \in \mathbb{Z}} \hat{\beta}_{lk} 2^{l/2} \psi(2^l y - k), \qquad y \in \mathbb{R},$$

where $\hat{\alpha}_{lk} = \int 2^{l/2} \phi(2^l x - k) \, dP_n(x)$, $\hat{\beta}_{lk} = \int 2^{l/2} \psi(2^l x - k) \, dP_n(x)$ and where $j_n \nearrow \infty$. This estimator was introduced in Doukhan and León (1990) and Kerkyacharian and Picard (1992). The latter authors proved—using wavelet theory as established by Daubechies (1992), Meyer (1992) and others—that this estimator is, for a suitable choice of $j_n$, an optimal estimator of $p_0$ in $\mathcal{L}^p$-loss, $1 \leq p < \infty$, if $p_0$ belongs to a Besov space $B_{pq}^t(\mathbb{R})$. Furthermore, "nonlinear" modifications of $p_n$ were shown to be optimal even in more general settings, including, in particular, the case when $t$ is unknown [see Donoho, Johnstone, Kerkyacharian and Picard (1995, 1996), Delyon and Juditsky (1996), Kerkyacharian, Picard and Tribouley (1996), Hall, Kerkyacharian and Picard (1998), Juditsky and Lambert-Lacroix (2004) and others]. The linear estimator is part of the analysis of these more complex nonlinear estimators. We refer to the monographs Härdle, Kerkyacharian, Picard and Tsybakov (1998) and Vidakovic (1999) for a general treatment of the use of wavelets in statistics.

In this article, we have three main goals: the first two consist in studying the limiting behavior of the linear estimator $p_n(y)$ both as an estimator for the true density function $p_0(y)$ and as an estimator $F_n^W(s) = \int_{-\infty}^s p_n(y) \, dy$ for the true distribution function $F(s) = \int_{-\infty}^s p_0(y) \, dy$, in sup-norm loss. Third—as a statistical application—we consider the same problems for a nonlinear modification of $p_n$, namely the "hard thresholding" wavelet density estimator.

In the first case, we show that under mild conditions,

(2)
$$\sup_{y \in \mathbb{R}} |p_n(y) - E p_n(y)| = O_{\text{a.s.}}\left(\sqrt{\frac{j_n 2^{j_n}}{n}}\right),$$

in fact we obtain an exact law of the logarithm for a suitably scaled version of $p_n - E p_n$, somewhat analogous to that of Deheuvels (2000) and Giné and Guillou (2002) for the Rosenblatt–Parzen kernel density estimator. A corollary of the proof also recovers, under weaker conditions, a result of Massiani (2003), where the supremum is taken over a bounded interval, as in the classical law of the logarithm of Stute (1982) for the Rosenblatt–Parzen estimator. The result (2) implies that, if $p_0$ is in the Besov space $B_{\infty\infty}^t(\mathbb{R})$ (or in the corresponding Hölder space of order $t$), then

(3)
$$\sup_{y \in \mathbb{R}} |p_n(y) - p_0(y)| = O_{\text{a.s.}}\left(\left(\frac{\log n}{n}\right)^{t/(2t+1)}\right),$$

if one chooses $2^{j_n} \simeq (n/\log n)^{1/(2t+1)}$, which is the optimal rate of convergence in sup-norm loss. These results are complemented by expectation bounds and convergence of Laplace transforms.

In the second case, we show, for $j_n$ as in the previous paragraph (and other choices), that the processes

$$\sqrt{n}(F_n^W - F)(s), \qquad s \in \mathbb{R}, \tag{4}$$

converge in law in the Banach space of bounded functions on $\mathbb{R}$ to the $P$-Brownian bridge process, and that

$$\sup_{s \in \mathbb{R}} |F_n^W(s) - F(s)| = O_{\text{a.s.}}\left(\sqrt{\frac{\log \log n}{n}}\right),$$

in fact, we obtain an exact law of the iterated logarithm and a strong approximation result. More generally, we then also prove uniform central limit theorems for the processes

$$\sqrt{n} \int_{\mathbb{R}} (p_n(y) - p_0(y)) f(y) \, dy, \qquad f \in \mathcal{F},$$

for several (Donsker) classes of functions $\mathcal{F}$. These results again parallel limit theorems for the classical Rosenblatt–Parzen estimator [see Bickel and Ritov (2003), Giné and Nickl (2008, 2009)].

To motivate the relevance of our third goal, note that the resolution $j_n$ under which the linear estimator achieves the optimal rate (3) for $p_0 \in B_{\infty\infty}^t(\mathbb{R})$ depends on $t$, which is typically unknown. To remedy this, Donoho et al. (1996) introduced (soft and hard) thresholding wavelet estimators: one first chooses $j_n$ sufficiently large and independent of $t$, and then deletes the wavelet coefficients $\hat{\beta}_{lk}$ in (1) in a certain range of $l$'s if they are smaller than a certain threshold. This estimator does not depend on $t$ anymore, but still achieves rates of convergence in the $\mathcal{L}^p$-loss, $1 \leq p < \infty$, that are optimal up to a logarithm factor, uniformly over compactly supported densities that are contained in balls of Besov spaces $B_{pq}^t(\mathbb{R})$, with $t$ unknown (but bounded). We show, as an application of our results for the linear estimator, that their hard thresholding estimator is *exact rate adaptive in the supnorm*, that is, it achieves the optimal rate (2) in the sup-norm, even without a logarithmic penalty, for (not necessarily compactly supported) $p_0$ in $B_{\infty\infty}^t(\mathbb{R})$, and any unspecified (but bounded) $t$. (In fact, this implies optimality over balls of densities in $B_{pq}^t(\mathbb{R})$ as well, cf. Remark 8 below.) Moreover, we prove that the hard thresholding wavelet density estimator also satisfies the central limit theorem (4). Hence this remarkable estimator is not only rate-adaptive in sup-norm loss, but also satisfies Bickel and Ritov's (2003) plug-in property.

The linear estimator in (1) can be expressed as a generalized kernel-type estimator

$$p_n(y) = \frac{2^j}{n} \sum_{i=1}^{n} K(2^j X_i, 2^j y),$$

where $K(x,y)$ is the wavelet projection kernel. It is interesting to compare to other kernel choices. The classical case would be the Parzen–Rosenblatt kernel density estimator, where $K(x,y) = K(x-y)$ with $K$ some probability density: if one makes the usual conversion from bandwidth $h$ to $2^{-j}$, one can compare directly with the classical kernel case, and we discuss this in some detail in Remark 6 below. In a nutshell, while the proof in the wavelet case follows a pattern similar to the one for classical kernels, some fundamental difficulties arise due to the fact that $K(x,y)$ is *not* a "convolution-type" kernel $K(x-y)$. Most importantly, the size of the random fluctuations of the (centered) wavelet estimator $p_n(y) - Ep_n(y)$ depends on $y$ not only through $p_0(y)$, but also through the quantity $\int K^2(2^j y, u)\, du$, which is part of the variance term, and which has periodic oscillations on $\mathbb{R}$ (unless one restricts oneself to the Haar wavelet). Among other things, this requires a normalization in the law of the logarithm that is not necessary in the convolution-kernel case of Stute (1982) and Giné and Guillou (2002). One might also be interested in considering projection kernels associated with other orthonormal systems that are not of wavelet type, as, for example, the Dirichlet kernel (which corresponds to an estimator based on Fourier series expansions). While our techniques may apply there as well, these kernels are often less interesting for estimating a function in the sup-norm, because of approximation-theoretic reasons: for example, the Fourier series of a uniformly continuous function might not converge at all points, and even if it does, its approximation properties in supnorm can be suboptimal.

Our proofs are based on techniques from empirical process theory. Note that if $p_0$ is compactly supported, or if $y$ varies in a fixed compact set, then $p_n(y) - Ep_n(y)$ consists of a *finite* sum of centered empirical wavelet coefficients, and in this case "finite dimensional" probabilistic methods are sufficient to analyze the limiting behavior of $p_n$ in the sup-norm. Otherwise, empirical process methods seem to be required. We show that the classes of functions naturally associated to wavelet density estimators are of Vapnik–Červonenkis type, and this allows the effective use of exponential inequalities for empirical processes [Talagrand (1996)] and entropy-based moment bounds [e.g., see Einmahl and Mason (2000), Giné and Guillou (2001)]. We also use that bounded subsets of Besov spaces are $P$-Donsker classes of functions [Nickl and Pötscher (2007)]. Wavelet theory is used throughout, and a brief summary of what we need precedes the main results.

## 2. Basic setup.

2.1. *Notation.* For an arbitrary (nonempty) set $M$, $\ell^\infty(M)$ will denote the Banach space of bounded real-valued functions $H$ on $M$ normed by $\|H\|_M := \sup_{m \in M} |H(m)|$, but we will use the symbol $\|h\|_\infty$ to denote $\sup_{x \in \mathbb{R}} |h(x)|$ for $h : \mathbb{R} \to \mathbb{R}$. For Borel-measurable functions $h : \mathbb{R} \to \mathbb{R}$ and Borel measures $\mu$ on $\mathbb{R}$, we set $\mu h := \int_{\mathbb{R}} h \, d\mu$, and we denote by $\mathcal{L}^p(\mu) := \mathcal{L}^p(\mathbb{R}, \mu)$ the usual Lebesgue-spaces of real-valued functions, normed by $\|\cdot\|_{p,\mu}$. If $d\mu(x) = dx$ is Lebesgue measure, we set $\mathcal{L}^p(\mathbb{R}) := \mathcal{L}^p(\mathbb{R}, \mu)$, and, for $1 \leq p < \infty$, we abbreviate the norm by $\|\cdot\|_p$. Similarly $\ell^p := \ell^p(\mathbb{Z})$, $1 \leq p \leq \infty$, are the usual sequence spaces, and we also denote their norm by $\|\cdot\|_p$ in a slight abuse of notation. All integrals are over the real line unless stated otherwise.

Let $X_1, \ldots, X_n$ be i.i.d. random variables with common law $P$ on $\mathbb{R}$, and denote by $P_n = n^{-1} \sum_{i=1}^{n} \delta_{X_i}$ the empirical measure. We assume throughout that the variables $X_i$ are the coordinate projections of $(\mathbb{R}^{\mathbb{N}}, \mathcal{B}^{\mathbb{N}}, P^{\mathbb{N}})$, and we set $\Pr := P^{\mathbb{N}}$. The empirical process indexed by $\mathcal{F} \subseteq \mathcal{L}^2(\mathbb{R}, P)$ is given by $f \mapsto \sqrt{n}(P_n - P)f, f \in \mathcal{F}$. Convergence in law of random elements in $\ell^\infty(\mathcal{F})$ is defined as, for example, in Dudley (1999) or de la Peña and Giné (1999), and will be denoted by the symbol $\rightsquigarrow_{\ell^\infty(\mathcal{F})}$. The class $\mathcal{F}$ is said to be $P$-*Donsker* if the centered Gaussian process $G_P$ with covariance $EG_P(f)G_P(g) = P[(f - Pf)(g - Pg)]$ is sample-bounded and sample-continuous w.r.t. the covariance semimetric, and if $\sqrt{n}(P_n - P) \rightsquigarrow_{\ell^\infty(\mathcal{F})} G_P$.

2.2. *Wavelet expansions and estimators.* We recall here some standard facts from wavelet theory [e.g., see Meyer (1992), Daubechies (1992), Härdle et al. (1998) or Vidakovic (1999)]. Let $\phi \in \mathcal{L}^2(\mathbb{R})$ be a father wavelet, that is, $\phi$ is such that $\{\phi(\cdot - k) : k \in \mathbb{Z}\}$ is an orthonormal system in $\mathcal{L}^2(\mathbb{R})$, and moreover the linear spaces $V_0 = \{f(x) = \sum_k c_k \phi(x - k) : \{c_k\}_{k \in \mathbb{Z}} \in \ell^2\}$, $V_1 = \{h(x) = f(2x) : f \in V_0\}, \ldots, V_j = \{h(x) = f(2^j x) : f \in V_0\}, \ldots$, are nested ($V_{j-1} \subseteq V_j$ for $j \in \mathbb{N}$) and such that $\bigcup_{j \geq 0} V_j$ is dense in $\mathcal{L}^2(\mathbb{R})$. For $\phi$ with compact support and

$$(5) \qquad K(y, x) := K_\phi(y, x) = \sum_{k \in \mathbb{Z}} \phi(y - k)\phi(x - k),$$

the functions $K_j(y, x) := 2^j K(2^j y, 2^j x), j \in \mathbb{N} \cup \{0\}$, are the kernels of the orthogonal projections of $\mathcal{L}^2(\mathbb{R})$ onto $V_j$, and we write $K_j(f)(y) = \int K_j(y, x) \times f(x) \, dx$ for this projection. We will use the following properties repeatedly throughout the proofs: if $\phi$ (not necessarily a father wavelet) is bounded and compactly supported, we have [e.g., Härdle et al. (1998), Lemma 8.6]

$$(6) \qquad |K(y, x)| \leq \Phi(y - x) \quad \text{and} \quad \sum_k |\phi(\cdot - k)| \in \mathcal{L}^\infty(\mathbb{R}),$$

where $\Phi : \mathbb{R} \to \mathbb{R}^+$ is bounded, compactly supported and symmetric. Furthermore, if $\phi$ is a bounded and compactly supported father wavelet, then, for every $x$,

$$(7) \qquad \int K(x, y) \, dy = 1$$

[see Corollary 8.1 in Härdle et al. (1998)]; moreover, for $f \in \mathcal{L}^p(\mathbb{R})$, $1 \le p \le \infty$, and fixed $j$, the series

$$K_j(f)(y) = \sum_{k \in \mathbb{Z}} 2^j \phi(2^j y - k) \int \phi(2^j x - k) f(x)\, dx, \qquad y \in \mathbb{R},$$

converges pointwise (since for each $y$ this is a finite sum). For $f \in \mathcal{L}^1(\mathbb{R})$, which is the main case in this article, the convergence of the series in fact takes place in $\mathcal{L}^p(\mathbb{R})$, $1 \le p \le \infty$. [For the reader's convenience, here is a proof: since $j$ is fixed, we can assume $j = 0$. Setting $a_k = \int \phi(x - k) f(x)\, dx$ we have $\int K_0(f)(x) \phi(x - k)\, dx = a_k$ by compactness of the support of $\phi$ and orthogonality, hence

$$\begin{aligned}\sum_k |a_k| &\le \int \sum_k |K_0(f)(x) \phi(x - k)|\, dx \le \sup_x \sum_k |\phi(x - k)| \|K_0(f)\|_1 \\ &\le c_1 \|\Phi * |f|\|_1 \le c_2 \|f\|_1\end{aligned} \tag{8}$$

by (6). Therefore, for any $1 \le p \le \infty$, $\sum_k \|a_k \phi(\cdot - k)\|_p \le \|\phi\|_p \sum_k |a_k| < \infty$.]

If now $\phi$ is a father wavelet and $\psi$ the associated mother wavelet so that $\{\phi(\cdot - k), 2^{l/2} \psi(2^l(\cdot) - k) : k \in \mathbb{Z}, l \in \mathbb{N} \cup \{0\}\}$ is an orthonormal basis of $\mathcal{L}^2(\mathbb{R})$ [see, e.g., Härdle et al. (1998), page 27], then any $f \in \mathcal{L}^p(\mathbb{R})$ has the formal expansion

$$f(y) = \sum_k \alpha_k(f) \phi(y - k) + \sum_{l=0}^{\infty} \sum_k \beta_{lk}(f) \psi_{lk}(y), \tag{9}$$

where $\psi_{lk}(y) = 2^{l/2} \psi(2^l y - k)$, $\alpha_k(f) = \int f(x) \phi(x - k)\, dx$, $\beta_{lk}(f) = \int f(x) \times \psi_{lk}(x)\, dx$. Since $(K_{l+1} - K_l) f = \sum_k \beta_{lk}(f) \psi_{lk}$ [e.g., Härdle et al. (1998), page 92], the partial sums of the series (9) are in fact given by

$$K_j(f)(y) = \sum_k \alpha_k(f) \phi(y - k) + \sum_{l=0}^{j-1} \sum_k \beta_{lk}(f) \psi_{lk}(y) \tag{10}$$

and—just as in the previous paragraph—one shows that, if $\phi, \psi$ are bounded and have compact support, then (10) converges pointwise and also in $\mathcal{L}^p(\mathbb{R})$, $1 \le p \le \infty$, if $f \in \mathcal{L}^1(\mathbb{R})$. If $p < \infty$, and $f \in \mathcal{L}^p(\mathbb{R})$, then convergence in (10) takes place in $\mathcal{L}^p(\mathbb{R})$ by a similar argument. Now using (6), (7), Minkowski's inequality for integrals and continuity of translations in $\mathcal{L}^p(\mathbb{R})$, we have $\|K_j(f) - f\|_p \le \int \Phi(u) \|f(2^{-j} u + \cdot) - f\|_p\, du \to 0$ as $j \to \infty$ for all $f \in \mathcal{L}^p(\mathbb{R})$, $1 \le p < \infty$, so that convergence of the wavelet series in (9) takes place in $\mathcal{L}^p(\mathbb{R})$.

Some regularity conditions on the wavelets $\phi, \psi$ will be needed. They parallel the order and moment conditions for convolution kernels in classical kernel density estimation. The standard conditions read as follows. Recall that $D\phi$ is the weak derivative of $\phi$ if $\int \phi Df = -\int (D\phi) f$ holds for all compactly supported infinitely differentiable functions $f : \mathbb{R} \to \mathbb{R}$.

CONDITION 1. (S). We say that the orthonormal system $\{\phi(\cdot - k), \psi_{lk} : k \in \mathbb{Z}, l \in \mathbb{N} \cup \{0\}\}$ is $S$-regular, if $\phi$ and $\psi$ are bounded and have compact support, and, if in addition, *one* of the following two conditions is satisfied: either (i) the father wavelet $\phi$ has weak derivatives up to order $S$ that are in $\mathcal{L}^p(\mathbb{R})$ for some $1 \leq p \leq \infty$; or (ii) the mother wavelet $\psi$ associated to $\phi$ satisfies $\int x^i \psi(x)\,dx = 0$, $i = 0, \ldots, S$.

The Haar wavelets, corresponding to $\phi = 1_{(0,1]}$ and $\psi = 1_{(0,1/2]} - 1_{(1/2,1]}$, satisfy this condition only for $S = 0$. And, for any given $S$ there exist compactly supported wavelets $\phi$ and $\psi$ that satisfy condition (S) [e.g., Daubechies' wavelets, see Daubechies (1992), Chapter 6, or Härdle et al. (1998)].

Given $X_1, \ldots, X_n$ i.i.d. with common absolutely continuous law $P$ on $\mathbb{R}$, the linear wavelet density estimator has the form

$$
\begin{aligned}
p_n(y) &:= P_n(K_{j_n}(y, \cdot)) = \frac{1}{n} \sum_{i=1}^n K_{j_n}(y, X_i) \\
&= \sum_k \hat{\alpha}_k \phi(y - k) + \sum_{l=0}^{j_n - 1} \sum_k \hat{\beta}_{lk} \psi_{lk}(y), \qquad y \in \mathbb{R},
\end{aligned}
$$
(11)

where $K$ is as in (5), $j_n \in \mathbb{N}$ satisfies $j_n \nearrow \infty$ as $n \to \infty$, and where $\hat{\alpha}_k = \int \phi(x - k)\,dP_n(x)$, $\hat{\beta}_{lk} = \int \psi_{lk}(x)\,dP_n(x)$ are the empirical wavelet coefficients. We note that for $\phi, \psi$ compactly supported, there are only finitely many $k$s for which these coefficients are nonzero (with the set of coefficients depending on $y$). Note that, if $\phi = 1_{(0,1]}$, then $p_n$ is just the usual histogram density estimator (with dyadic binpoints). For general compactly supported wavelets $\phi, \psi$, the estimator $p_n$ was first studied by Doukhan and León (1990) and Kerkyacharian and Picard (1992).

2.3. *Besov spaces.* To deal with the approximation error ("bias term") of wavelet density estimators, and for some proofs, we introduce the Besov spaces $B_{pq}^s(\mathbb{R})$, which form a general scale of smooth function spaces (that contain all the classical ones as special cases). Besov spaces can be defined in several equivalent ways, the classical one being in terms of $\mathcal{L}^p$–$\mathcal{L}^q$ norms of the second differences $|h|^{-sq-1} \times (D^{s-\{s\}} f(\cdot + h) + D^{s-\{s\}} f(\cdot - h) - 2D^{s-\{s\}} f(\cdot))$ of weak derivatives of $f$, where $0 < \{s\} \leq 1$ and $s - \{s\} \in \mathbb{N} \cup \{0\}$. Wavelet bases provide another characterization of these spaces, hence it is most convenient for our purposes to define them in this way.

DEFINITION 1. Let $1 \leq p, q \leq \infty$, $0 < s < S$, $s \in \mathbb{R}$, $S \in \mathbb{N}$. Let $\phi$ be a bounded, compactly supported father wavelet that satisfies part (i) of Condition 1(S), and denote by $\alpha_k(f)$ and $\beta_{lk}(f)$ the wavelet coefficients of $f \in \mathcal{L}^p(\mathbb{R})$. The Besov space $B^s_{pq}(\mathbb{R})$ is defined as the set of functions

$$\left\{ f \in \mathcal{L}^p(\mathbb{R}) : \|f\|_{s,p,q} := \|\alpha_{(\cdot)}(f)\|_p \right.$$
$$\left. + \left( \sum_{l=0}^{\infty} (2^{l(s+1/2-1/p)} \|\beta_{l(\cdot)}(f)\|_p)^q \right)^{1/q} < \infty \right\}$$

with the obvious modification in case $q = \infty$.

REMARK 1 (Properties of Besov spaces). That this definition coincides with the more classical ones follows, for instance, from Meyer (1992, page 200) or Section 9 in Härdle et al. (1998). The definition is independent of the choice of $\phi, \psi$, and one has the continuous imbedding of $B^r_{pq}(\mathbb{R})$ [defined via $\phi$ satisfying part (i) of Condition 1(R) with $0 < r < R$] into $B^s_{pq}(\mathbb{R})$ [defined via a possibly different $\phi'$ satisfying part (i) of Condition 1(S) with $0 < s < S$ for $r \geq s$]. We also recall some well-known relations of $B^s_{pq}(\mathbb{R})$ to classical smooth function spaces [see, e.g., Triebel (1983)]: for example, $B^s_{pq}(\mathbb{R})$ is continuously imbedded into $\mathcal{L}^p(\mathbb{R})$ for $1 \leq p \leq \infty$, and, if $\mathsf{C}^s(\mathbb{R})$ are the classical Hölder spaces (of $s$-times continuously differentiable functions in case $s \in \mathbb{N}$), then

(12) $$B^s_{\infty 1}(\mathbb{R}) \hookrightarrow \mathsf{C}^s(\mathbb{R}) \hookrightarrow B^s_{\infty \infty}(\mathbb{R})$$

holds, where the second imbedding is even an identity if $s$ is noninteger; and one also has the Sobolev type imbedding $B^s_{pq}(\mathbb{R}) \hookrightarrow \mathsf{C}^{s-1/p}(\mathbb{R})$ for $s > 1/p$ or $s = 1/p$ and $q = 1$. Further examples are the classical Sobolev spaces $H^s(\mathbb{R}) = \{f \in \mathcal{L}^2(\mathbb{R}) : |Ff(\cdot)|^2(1 + |\cdot|^2)^s \in \mathcal{L}^2(\mathbb{R})\}$, where $F$ is the Fourier transform, for which one has $H^s(\mathbb{R}) = B^s_{22}(\mathbb{R})$; and if $BV(\mathbb{R}) = \{f : v_1(f) < \infty\}$, where $v_1$ is defined in (13) below, then $B^1_{11}(\mathbb{R}) \hookrightarrow BV(\mathbb{R}) \cap \mathcal{L}^1(\mathbb{R}) \hookrightarrow B^1_{1\infty}(\mathbb{R})$.

**3. Entropy and expectation bounds.** In this section we will show that certain classes of functions related to the kernel $K(y, x) = \sum_{k \in \mathbb{Z}} \phi(y-k)\phi(x-k)$ are VC-type classes of functions, meaning that they have $\mathcal{L}^2(Q)$ covering numbers of polynomial order, uniformly in all probability measures $Q$. Using expectation inequalities for VC-classes, we obtain—as an immediate consequence—a finite sample inequality for the expected value of the deviation of the wavelet estimator from its mean. Also, these VC-bounds will be applied in later sections to obtain various exponential inequalities for wavelet density estimators.

A function $h$ is of bounded $p$-variation on $\mathbb{R}$, $0 < p < \infty$, if

(13)
$$v_p(h) := \sup\left\{\sum_{i=1}^n |f(x_i) - f(x_{i-1})|^p : \right.$$
$$\left. n \in \mathbb{N}, -\infty < x_0 < x_1 < \cdots < x_n < \infty\right\}$$

is finite. The following lemma—which uses (and generalizes) a result due to Nolan and Pollard (1987)—will be useful in what follows. As usual, for $\mathcal{H}$ a class of functions in $\mathcal{L}^r(Q)$, $1 \leq r < \infty$, $N(\mathcal{H}, \mathcal{L}^r(Q), \varepsilon)$ denotes the minimal number of $\mathcal{L}^r(Q)$-balls of radius less than or equal to $\varepsilon$, that cover $\mathcal{H}$. The logarithm of the covering number is the $\mathcal{L}^r(Q)$-metric entropy of $\mathcal{H}$.

LEMMA 1. *Let $h : \mathbb{R} \to \mathbb{R}$ be a function of bounded $p$-variation, $p \geq 1$. Set*
$$\mathcal{H} = \{h((\cdot)t - s) : t, s \in \mathbb{R}\}.$$
*Then $\mathcal{H}$ satisfies*
$$\sup_Q N(\mathcal{H}, \mathcal{L}^2(Q), \varepsilon) \leq \left(\frac{A}{\varepsilon}\right)^v, \qquad 0 < \varepsilon < A,$$
*with finite positive constants $A, v$ depending only on $h$, and where the supremum extends over all Borel probability measures $Q$ on $\mathbb{R}$.*

PROOF. It is known that $h$ is equal to $g \circ f$ where $f$ is nondecreasing with range contained in $[0, v_p(h)]$ and $g$ is $1/p$-Hölder-continuous on the full interval $[0, v_p(h)]$ [see Love and Young (1937) and also Dudley (1992), page 1971]. The set $\mathcal{M}$ of dilations and translations of $f$ satisfies the required entropy bound with $\mathcal{L}^2(Q)$ replaced by $\mathcal{L}^r(Q)$ for any $r > 0$ (where $\|\cdot\|_{r,Q} = \int |\cdot|^r \, dQ$ if $r < 1$), with a constant $A$ that depends only on $r$ times $v_1(f)$ [see Nolan and Pollard (1987) and de la Peña and Giné (1999), page 224, for $r < 1$]. Since
$$\int |g(m_1) - g(m_2)|^2 \, dQ \leq \int |m_1 - m_2|^{2/p} \, dQ,$$
it follows that any $\varepsilon$-covering of $\mathcal{M}$ for $\mathcal{L}^{2/p}(Q)$ induces a $\varepsilon^s$-covering of $\mathcal{H}$ of the same cardinality, for $s = 1/p$ if $2/p \geq 1$ and $s = 1/2$ otherwise, proving the lemma (for suitable $v$ depending only on $p$). □

We will impose the following condition on the function $\phi$ defining the kernel $K$ in (5).

CONDITION 2. $\phi : \mathbb{R} \to \mathbb{R}$ is of bounded $p$-variation for some $1 \leq p < \infty$ and vanishes on $(B_1, B_2]^c$ for some $-\infty < B_1 < B_2 < \infty$.

The Haar father wavelet $\phi = 1_{(0,1]}$ is of bounded variation ($p=1$) and hence satisfies Condition 2. Furthermore, since any $\alpha$-Hölder-continuous function ($0 < \alpha \leq 1$) with compact support is also of bounded $1/\alpha$-variation, Condition 2 is also satisfied, for example, for all Daubechies' (father) wavelets [see, e.g., Härdle et al. (1998), Remark 7.1]. It should be noted that not all Daubechies' wavelets are of bounded variation, but they are all Hölder continuous for some $0 < \alpha < 1$, which is why the generalization to $p$-variation of the result of Nolan and Pollard (1987) is useful in the present context.

Now for $\phi$ satisfying Condition 2, define

$$\mathcal{F}_\phi = \left\{ \sum_{k \in \mathbb{Z}} \phi(2^j y - k) \phi(2^j(\cdot) - k) : y \in \mathbb{R}, j \in \mathbb{N} \cup \{0\} \right\} \tag{14}$$

and

$$\mathcal{D}_{\phi,j} = \left\{ \sum_{k \in \mathbb{Z}} 2^j \int_{-\infty}^t \phi(2^j y - k) \, dy \, \phi(2^j(\cdot) - k) : t \in \mathbb{R} \right\}, \qquad j \in \mathbb{N} \cup \{0\}. \tag{15}$$

Notice that by (6), both classes have a constant envelope.

LEMMA 2. *Let $\mathcal{G}$ be either $\mathcal{F}_\phi$ or $\mathcal{D}_{\phi,j}$, where $\phi$ satisfies Condition 2. Then we have the uniform entropy bound*

$$\sup_Q N(\mathcal{G}, \mathcal{L}^2(Q), \varepsilon) \leq \left( \frac{A}{\varepsilon} \right)^v, \qquad 0 < \varepsilon < A \tag{16}$$

*for $A, v$ positive and finite constants depending only on $\phi$ (and not on $j$ for $\mathcal{D}_{\phi,j}$), and where the supremum extends over all Borel probability measures $Q$ on $\mathbb{R}$.*

PROOF. The case of $\mathcal{F}_\phi$: for $y, j$ fixed, the sum $\sum_{k \in \mathbb{Z}} \phi(2^j y - k) \phi(2^j(\cdot) - k)$ consists of at most $[B_2 - B_1] + 1$ summands, each of which has the form

$$\phi(2^j y - k) \phi(2^j(\cdot) - k) = c_{j,y,k} \phi(2^j(\cdot) - k),$$

where $k$ is a fixed integer satisfying $2^j y - B_2 \leq k < 2^j y - B_1$, and where $|c_{j,y,k}| \leq \|\phi\|_\infty$. Since $\phi$ is of bounded $p$-variation, Lemma 1 above applies to the class $\mathcal{M}$ of dilations and translations of $\phi$, yielding the entropy bound (16) for $\mathcal{M}$ (with different constants $A, v$). The class $\mathcal{F}_\phi$ consists of linear combinations of at most $[B_2 - B_1] + 1$ elements of $\mathcal{M}$, whose coefficients are bounded in absolute value by $\|\phi\|_\infty$. For given $\varepsilon' > 0$, take an $\varepsilon'$-dense subset $\{a_l\}$ of $[-\|\phi\|_\infty, \|\phi\|_\infty]$ and an $L^2(Q)$-$\varepsilon'$-dense subset $\{m_i(\cdot)\}$ of $\mathcal{M}$. Then $\{\sum_{k=1}^{[B_2-B_1]+1} a_{l_k} m_{i_k}(\cdot)\}_{l,i}$ are the centers of a covering of $\mathcal{F}_\phi$ by $L_2(Q)$ balls of radius $\varepsilon = ([B_2 - B_1] + 1)(\|\phi\|_\infty + 1)\varepsilon'$, and a simple computation on covering numbers shows that the required entropy bound holds for $\mathcal{F}_\phi$.

The case of $\mathcal{D}_{\phi,j}$: by the support assumption on $\phi$, we have for every fixed $t$,

$$\sum_{k \in \mathbb{Z}} 2^j \int_{-\infty}^t \phi(2^j y - k)\, dy\, \phi(2^j(\cdot) - k)$$
$$= c \sum_{k \leq 2^j t - B_2} \phi(2^j(\cdot) - k) + \sum_{2^j t - B_2 < k \leq 2^j t - B_1} c_{j,t,k} \phi(2^j(\cdot) - k),$$

where $c = \int_{-\infty}^\infty \phi(y)\, dy$ and $|c_{j,t,k}| \leq \|\phi\|_1$. The class of functions

$$\left\{ \sum_{2^j t - B_2 < k < 2^j t - B_1} c_{j,t,k} \phi(2^j(\cdot) - k) : t \in \mathbb{R} \right\}$$

satisfies the bound (16) with $A, v$ independent of $j$, by the argument in the first part of the proof. Each function in the class

$$\left\{ c \sum_{k \leq 2^j t - B_2} \phi(2^j(\cdot) - k) : t \in \mathbb{R} \right\}$$

is the difference of two functions, one in each of the classes

$$\left\{ c \sum_{k \leq 2^j t - B_2} \phi_+(2^j(\cdot) - k) : t \in \mathbb{R} \right\}, \quad \left\{ c \sum_{k \leq 2^j t - B_2} \phi_-(2^j(\cdot) - k) : t \in \mathbb{R} \right\},$$

where $\phi = \phi_+ - \phi_-$ and $\phi_+, \phi_- \geq 0$. But these classes are linearly ordered, so their subgraphs are ordered by inclusion, and therefore are VC-subgraph of index 1 [cf. Dudley (1999), Theorem 4.2.6.]. The entropy bound for $D_{\phi,j}$ follows from these observations and, again, a simple computation on covering numbers. $\square$

Using expectation bounds for VC-classes of functions [e.g., Einmahl and Mason (2000), Giné and Guillou (2001)], the last lemma already implies the following result.

PROPOSITION 1. *Let $K(y,x) = \sum_k \phi(y-k)\phi(x-k)$ where $\phi$ satisfies Condition 2. Suppose that $P$ has a bounded density $p_0$. Let $0 < L < \infty$. Then, for every $n \in \mathbb{N}$ and every fixed integer $l \geq 0$, $l' = \max(l,1)$, there exists a fixed constant $c$ independent of $n, l$ such that*

(17) $$\sup_{p_0 : \|p_0\|_\infty \leq L} E \sup_{y \in \mathbb{R}} |(P_n - P)K_l(y, \cdot)| \leq c(\sqrt{2^l l'/n} + (2^l l'/n)).$$

*If, in addition, $\phi$ is a father and $\psi$ an associated mother wavelet satisfying Condition 1(0), then, setting $\hat{\alpha}_k = \int \phi(x-k) \, dP_n(x)$, $\hat{\beta}_{lk} = \int \psi_{lk}(x) \, dP_n(x)$, $\alpha_k = \int \phi(x-k) \, dP(x)$ and $\beta_{lk} = \int \psi_{lk}(x) \, dP(x)$, we have*

(18)
$$\sup_{p_0 : \|p_0\|_\infty \leq L} E \sup_{k \in \mathbb{Z}} |\hat{a}_k - a_k| \leq C/\sqrt{n},$$

$$\sup_{p_0 : \|p_0\|_\infty \leq L} E \sup_{k \in \mathbb{Z}} |\hat{\beta}_{lk} - \beta_{lk}| \leq C(\sqrt{l'/n} + 2^{l/2} l'/n)$$

*for all $l \geq 0$ and a constant $C$ independent of $l, n$.*

PROOF. We apply Lemma 2 and the expectation inequality (57) below to the class

$$\mathcal{F} = \{K(2^l y, 2^l(\cdot)) - P(K(2^l y, 2^l(\cdot))) : y \in \mathbb{R}\},$$

which satisfies the same entropy bound as $\mathcal{F}_\phi$, and which has constant envelope $U$ independent of $p_0$ [using (6)]. To bound second moments, we use (6) to the effect that

$$\sup_y \int K^2(2^l y, 2^l x) p_0(x) \, dx = \sup_y 2^{-l} \int K^2(2^l y, 2^l y + u) p_0(y + 2^{-l} u) \, du$$

$$\leq 2^{-l} \|p_0\|_\infty \sup_{y \in \mathbb{R}} \|K(y, y + \cdot)\|_2^2 \leq 2^{-l} L \|\Phi\|_2^2 = \sigma^2.$$

The first claim of the proposition now follows from (57) [and the measurability Remark 2 below, which implies that the supremum over $y$ in (17) is in fact countable]. For the second claim, set $\tilde{\beta}_{lk} = \hat{\beta}_{lk} - \beta_{lk}$, define $h(x) := \sum_r \tilde{\beta}_{lr} \psi_{lr}(x) = (P_n - P)(K_{l+1} - K_l)(x)$, and note that

$$\int h \psi_{lk} = \int \sum_r \tilde{\beta}_{lr} \psi_{lr} \psi_{lk} = \tilde{\beta}_{lk}$$

(since the sum has finitely many terms and the $\psi_{lr}$ are orthogonal). Consequently,

$$\sup_k |\tilde{\beta}_{lk}| \leq \|h\|_\infty \sup_k 2^{l/2} \int |\psi(2^l x - k)| \, dx$$

$$\leq 2^{-l/2} \|\psi\|_1 \|(P_n - P)(K_{l+1} - K_l)\|_\infty,$$

which gives the bound (18) for $\hat{\beta}_{lk} - \beta_{lk}$ by the first part. In the case of $\hat{\alpha}_k - \alpha_k$ we have by a similar argument that $\sup_k |\hat{\alpha}_k - \alpha_k| \leq \|(P_n - P)K_0\|_\infty \|\phi\|_1$, and the result follows from the case $l = 0$ in the first part of the proposition. □

Inequalities analogous to (17) and (18) hold as well if the first moment is replaced by $p$th moments. This can proved either directly, or by combining the above moment bounds with Proposition 3.1 in Giné, Latala and Zinn (2000).

If $p_0$ has compact support, then the number of nonzero wavelet coefficients $\hat{\beta}_{lk}, \beta_{lk}(p_0)$ at level $l$ is finite and of the order $2^l$. In this case, the above proposition follows from Bernstein's inequality combined with a simple convexity argument that elaborates on one due to Pisier [cf. van der Vaart and Wellner (1996), Lemma 2.2.10]. However, if $p_0$ does not have compact support, empirical process methods seem to be unavoidable to prove (17) and (18).

Using the second part of Lemma 2, one can obtain a similar expectation bound for the distribution function $F_n^W(s) = \int_{-\infty}^s p_n(y)\,dy$ of the wavelet density estimator, as we do after Lemma 4 below.

**4. Rates of almost sure uniform consistency for the wavelet density estimator.** We will now derive best possible almost sure rates of convergence for the deviation of the estimator $p_n(y) = P_n(K_{j_n}(y,\cdot))$ from its mean $Ep_n(y) = P(K_{j_n}(y,\cdot)) = K_{j_n}(p_0)(y)$ uniformly in $y \in \mathbb{R}$. We also obtain a uniform law of the logarithm for a suitably scaled version of $p_n - Ep_n$. The results from this section are compared to similar results for the classical convolution kernel estimator in Remark 6 below.

For $K(y,x) = \sum_k \phi(y-k)\phi(x-k)$ as in (5), define the function

$$\bar{K}(y,x) := \frac{K(y,x)}{\sqrt{\int K^2(y,u)\,du}} = \frac{K(y,x)}{\sqrt{\sum_k \phi^2(y-k)}}. \tag{19}$$

Using (6) and (7), it is easy to see that, if $\phi$ is bounded and compactly supported then there exist finite non-zero constants $D_1, D_2$ independent of $y$ such that

$$D_1^2 \le \int K^2(y,u)\,du \le D_2^2. \tag{20}$$

We now proceed to prove the first main result, which is not the most exact, but requires minimal hypotheses. Let $j_n \nearrow \infty$ be a sequence of nonnegative integers satisfying the following conditions:

$$\frac{n}{j_n 2^{j_n}} \to \infty, \qquad \frac{j_n}{\log \log n} \to \infty, \qquad \sup_{n \ge n_0}(j_{2n} - j_n) \le \tau \tag{21}$$

for some $\tau \ge 1$ and some $n_0 < \infty$.

THEOREM 1. *Let $\phi$ be a father wavelet satisfying Condition 2. Suppose that $P$ has a bounded density $p_0$ and that $j_n$ satisfies (21). Then we have*

$$\limsup_{n \to \infty} \sqrt{\frac{n}{j_n 2^{j_n}}} \sup_{y \in \mathbb{R}} \left| \frac{p_n(y) - Ep_n(y)}{\sqrt{\sum_k \phi^2(2^{j_n} y - k)}} \right| = C \qquad a.s., \tag{22}$$

where $C^2 \leq M^2 2^\tau \|p_0\|_\infty$ for a constant $M$ that depends only on $D_1, D_2, \|\Phi\|_\infty$ [cf. (6)] and on the VC-characteristics $A$ and $v$ of $\mathcal{F}_\phi$.

PROOF. Let $n_k = 2^k$. We have

$$\Pr\left\{\max_{n_{k-1}<n\leq n_k} \sup_{y\in\mathbb{R}} \sqrt{\frac{1}{n2^{-j_n}j_n}} \left|\sum_{i=1}^n (\bar{K}(2^{j_n}y, 2^{j_n}X_i) - E\bar{K}(2^{j_n}y, 2^{j_n}X))\right| > s\right\}$$

$$(23) \quad \leq \Pr\left\{\max_{n_{k-1}<n\leq n_k} \sup_{\substack{y\in\mathbb{R}\\ j_{n_{k-1}}<j\leq j_{n_k}}} \left|\sum_{i=1}^n (\bar{K}(2^j y, 2^j X_i) - E\bar{K}(2^j y, 2^j X))\right| > s\sqrt{\frac{n_{k-1}j_{n_k}}{2^{j_{n_k}}}}\right\},$$

where $j \in \mathbb{N}$. To estimate the last probability, we will apply Talagrand's inequality to the classes of functions

$$\mathcal{F}_k = \{\bar{K}(2^j y, 2^j(\cdot)) - P(\bar{K}(2^j y, 2^j(\cdot))) : y \in \mathbb{R}, j_{n_{k-1}} < j_n \leq j_{n_k}\},$$

which have constant envelope $2\|\Phi\|_\infty/D_1$ [by (6) and (20)] and satisfy the same entropy bound as $\mathcal{F}_\phi$ in Lemma 2, with a possibly different $A, v$—but independent of $k$—by that lemma and a simple computation on covering numbers (since $(\int K^2(2^j y, u)\,du)^{-1} \in [1/D_2^2, 1/D_1^2]$ for all $y \in \mathbb{R}$). Consequently we may apply inequality (60) below with $U := U_k = 2\|\Phi\|_\infty/D_1$, and $\sigma^2 = 2^{-j_{n_{k-1}}}\|p_0\|_\infty$, where the bound on $\sigma$ follows from

$$(24) \quad \sup_y \int \bar{K}^2(2^j y, 2^j x)p_0(x)\,dx = \sup_y 2^{-j}\frac{\int K^2(2^j y, u)p_0(2^{-j}u)\,du}{\int K^2(2^j y, u)\,du}$$
$$\leq 2^{-j}\|p_0\|_\infty.$$

To be precise, we also need that the supremum in (23) is countable, and we show in Remark 2 below that this is the case. Setting $s = \sqrt{2^{\tau+2}C_1}\|p_0\|_\infty^{1/2}$ makes $t = s\sqrt{n_{k-1}2^{-j_{n_k}}j_{n_k}}$ an admissible choice in (60) for all $k$ large enough by the first and third conditions in (21). As a consequence, for these values of $k$, the probability in question is bounded from above by

$$R\exp\left\{-\frac{1}{C_3}\frac{s^2 n_{k-1}2^{-j_{n_k}}j_{n_k}}{n_k 2^{-j_{n_{k-1}}}\|p_0\|_\infty}\right\} \leq R\exp\left\{-\frac{2C_1 j_{n_k}}{C_3}\right\}.$$

Now the second limit in (21) becomes $j_{n_k}/\log k \to \infty$, hence the last expression is the general term of a convergent series. Thus, modulo measurability, we have proved that (for the stipulated $s$),

$$\sum_k \Pr\left\{\max_{n_{k-1}<n<n_k} \sup_{y\in\mathbb{R}} \sqrt{\frac{2^{j_n}}{nj_n}}\left|\sum_{i=1}^n (\bar{K}(2^{j_n}y, 2^{j_n}X_i) - E\bar{K}(2^{j_n}y, 2^{j_n}X))\right| > s\right\} < \infty,$$

which gives the theorem by Borel–Cantelli and the 0–1 law. □

REMARK 2 (Measurability). In order to apply Talagrand's inequality in the previous theorem, we must show that the supremum in (22) is in fact a countable supremum. Let $T_1$ be the set of discontinuities of $\phi$, which is countable since, $\phi$ being of bounded $p$-variation, it is the composition of a Hölder-continuous with a nondecreasing function (see the proof of Lemma 1). Let $T_0$ be a countable dense subset of $\mathbb{R} \setminus T_1$ and define $T = \{2^{-j}(z+k) : k \in \mathbb{Z}, z \in T_0 \cup T_1\}$. For each $y$, let $\phi_y = (\phi(2^j y - k) : k \in \mathbb{Z}) \in \ell^\infty(\mathbb{Z})$. We first prove that

$$\{\phi_y; y \in \mathbb{R}\} \subseteq \overline{\{\phi_y; y \in T\}}, \tag{25}$$

where the closure is in $\ell^\infty(\mathbb{Z})$. Given $y \in \mathbb{R}$, two cases are possible. Either $2^j y - k$ is a discontinuity point of $\phi$ for some $k \in \mathbb{Z}$, or $2^j y - k \in T_1^c$ for all $k$. In the first case, $y \in T$. In the second case, $\phi_y$ can be approximated by $\phi_{y_\delta}$, $y_\delta \in T$ as follows. Let $k_0$ be the largest integer such that $2^j y - k_0 > B_2$, and $k_N = k_0 + N$ be the smallest integer such that $2^j y - k_N < B_1$, and set $k_i = k_0 + i$, $i = 0, \ldots, N$. Note that $N \le B_2 - B_1 + 1$. Let $0 < \delta_0 < 1$ be such that $\phi$ is continuous on the neighborhood $N_i(\delta_0)$ of $2^j y - k_i$ of radius $\delta_0$, $i = 0, \ldots, N$. For $\delta \le \delta_0$ let $z \in N_0(\delta) \cap T_0$ and define $y_\delta = 2^{-j}(z + k_0) \in T$. Then $|2^j y - k_i - (2^j y_\delta - k_i)| < \delta$, $i = 0, \ldots, N$. Hence, by continuity of $\phi$ at $2^j y - k_i$, we have $\max_{0 \le i \le N} |\phi(2^j y - k_i) - \phi(2^j y_\delta - k_i)| \to 0$ as $\delta \to 0$. Since, moreover, $\phi(2^j y - k) = \phi(2^j y_\delta - k) = 0$ if $k \notin \{k_0, \ldots, k_N\}$, we have $\phi_{y_\delta} \to \phi_y$ in $\ell^\infty(\mathbb{Z})$ as $\delta \to 0$ concluding the proof of (25). Now

$$\sum_{i=1}^n (\bar{K}(2^j y, 2^j X_i) - E\bar{K}(2^j y, 2^j X)) = \frac{\sum_k \phi(2^j y - k) c_k}{\sqrt{\sum_k \phi^2(2^j y - k)}} =: \Gamma(y),$$

where $c_k$ are random variables satisfying $\sum_k |c_k| \le c < \infty$ for $c$ nonrandom by (6), and where $\sum_k \phi^2(2^j y - k) \ge D_1^2 > 0$ by (20). Hence if $\phi_{y_\delta} \to \phi_y$ in $\ell^\infty(\mathbb{Z})$ then $\Gamma(y_\delta) \to \Gamma(y)$ (as $\delta \to 0$). This, together with (25) proves that $\sup_{y \in \mathbb{R}} |\Gamma(y)| = \sup_{y \in T} |\Gamma(y)|$. That is, the supremum in (22) is countable.

REMARK 3. The proof of Theorem 1 also shows that, under the conditions of this theorem,

$$\limsup_{n \to \infty} \sqrt{\frac{n}{j_n 2^{j_n}}} \sup_{y \in \mathbb{R}} |p_n(y) - E p_n(y)| = C \quad \text{a.s.,} \tag{26}$$

where $C^2 \le M 2^\tau \|\Phi\|_2^2 \|p_0\|_\infty$. [The only difference is that in this case we use the variance estimate $\sigma^2 = 2^{-j} \|p_0\|_\infty \|\Phi\|_2^2$, which follows as in (24).]

The following corollary to the proof of Theorem 1 will be needed for the more exact result below.

COROLLARY 1. *Let $D \subset \mathbb{R}$ be such that $\|p_0\|_D > 0$. If, in addition to the hypotheses in Theorem 1, $p_0$ is uniformly continuous, then*

$$\limsup_{n\to\infty} \sqrt{\frac{n}{j_n 2^{j_n}}} \sup_{y\in D} \left| \frac{p_n(y) - Ep_n(y)}{\sqrt{\sum_k \phi^2(2^{j_n}y - k)}} \right| = C \quad \text{a.s.,}$$

*where $C^2 \leq M^2 2^\tau \|p_0\|_D$ and where $M$ is as in Theorem 1.*

PROOF. The proof is, as in Theorem 1, after observing that for every $\varepsilon > 0$ and $k$ large enough, the bound in (24), for $y \in D$, becomes $\sigma^2 = (1+\varepsilon)2^{-j_{n_k-1}}\|p_0\|_D$, by uniform continuity of $p_0$ and since, for all $x$, $K(x, x+u) = 0$ if $|u| > B_2 - B_1$. □

To obtain the exact constant in the almost sure limit, we now proceed to give a lower bound.

PROPOSITION 2. *Let $\phi$ be a bounded father wavelet vanishing on $(B_1, B_2)^c$, $-\infty < B_1 < B_2 < \infty$, and assume that $P$ has a bounded continuous density $p_0$. Then, if $j_n/\log\log n \to \infty$, we have*

$$\liminf_{n\to\infty} \sqrt{\frac{n}{(2\log 2)j_n 2^{j_n}}} \sup_{y\in\mathbb{R}} \left| \frac{p_n(y) - Ep_n(y)}{\sqrt{\sum_k \phi^2(2^{j_n}y - k)}} \right| \geq \|p_0\|_\infty^{1/2} \quad \text{a.s.}$$

PROOF. By Proposition 2 in Einmahl and Mason (2000), the conclusion holds if, for every $\varepsilon > 0$ and $n$ large enough (depending on $\varepsilon$), we can find $k_n = k_n(\varepsilon)$ points $z_{in} = z_{in}(\varepsilon)$, $i = 1, \ldots, k_n$, such that, if

$$g_i^{(n)}(x) = g_i^{(n,\varepsilon)}(x) = \bar{K}(2^{j_n}z_{in}, 2^{j_n}x),$$

then the following conditions hold (for all $n$ large enough and for constants $r$, $\mu_i$, $\sigma_i$, $i = 1, 2$, depending on $\varepsilon$):

(a) $\Pr\{g_i^{(n)}(X) \neq 0, g_k^{(n)}(X) \neq 0\} = 0$, $i \neq k$,
(b) $\sum_{i=1}^{k_n} \Pr\{g_i^{(n)}(X) \neq 0\} \leq 1/2$,
(c) $2^{-j_n}k_n \to r \in (0, \infty)$,
(d) $2^{-j_n}\mu_1 \leq E(g_i^{(n)}(X)) \leq 2^{-j_n}\mu_2$ for some $-\infty < \mu_1 < \mu_2 < \infty$,
(e) $\sigma_1 2^{-j_n/2} \leq \sqrt{\text{Var}(g_i^{(n)})(X)} \leq \sigma_2 2^{-j_n/2}$ for some $0 < \sigma_1 < \sigma_2 < \infty$,
(f) $\sup_{i,n} \|g_i^{(n)}\|_\infty < \infty$,
(g) $\lim_{\varepsilon\to 0} \sigma_1(\varepsilon) = \lim_{\varepsilon\to 0} \sigma_2(\varepsilon) = \|p_0\|_\infty$.

We proceed to verify these conditions. Given $\varepsilon > 0$, let $I$ be an interval such that $p_0(x) \geq (1-\varepsilon)\|p_0\|_\infty$ for all $x \in I$; and such that $\Pr\{X \in I\} \leq 1/2$. Such an interval exists because $p_0$ is bounded and continuous. Set $I = [a, b]$ and define

$$z_{in} = a + 3(B_2 - B_1)i2^{-j_n},$$

where
$$i = 1, 2, \ldots, \qquad \left[\frac{(b-a)2^{j_n}}{3(B_2 - B_1)}\right] - 1 := k_n.$$

For (a) note that $K(2^{j_n} z_{in}, 2^{j_n} x) \neq 0$ implies $|x - z_{in}| 2^{j_n} \leq B_2 - B_1$, and that $|z_{in} - z_{kn}| > 2^{-j_n} 3|B_2 - B_1|$ by construction, which together imply that the set in question is empty. For (b) note that by (a) the sum of the probabilities in (b) is $\Pr(\bigcup_{i=1}^{k_n} \{g_i^{(n)}(X) \neq 0\}) \leq \Pr\{X \in I\} \leq 1/2$. By construction, the limit in (c) is $\frac{b-a}{3(B_2-B_1)}$. Condition (f) follows immediately from (6) and the assumption on $\phi$. Conditions (d) and (e) are implied by the following estimates. First,

$$
\begin{aligned}
(27) \quad \int |\bar{K}(2^{j_n} z_{in}, 2^{j_n} x)| p_0(x)\, dx \\
\leq D_1^{-1} \int |K(2^{j_n} z_{in}, 2^{j_n} x)| p_0(x)\, dx \\
\leq 2^{-j_n} D_1^{-1} \int |K(2^{j_n} z_{in}, 2^{j_n} z_{in} + u)| p_0(z_{in} + u 2^{-j_n})\, du \\
\leq 2^{-j_n} D_1^{-1} \|p_0\|_\infty \|\Phi\|_1,
\end{aligned}
$$

where we use (20) in the first inequality and (6) in the last, and

$$\int \bar{K}^2(2^{j_n} z_{in}, 2^{j_n} x) p_0(x)\, dx \leq 2^{-j_n} \|p_0\|_\infty$$

by (24), which give the upper bounds in (d) and (e) with $\mu_2 = D_1^{-1} \|p_0\|_\infty \|\Phi\|_1$ and $\sigma_2^2 = \|p_0\|_\infty$. Second, for the lower bound in (d), again using (6), (7) and (20),

$$
\begin{aligned}
\int \bar{K}(2^{j_n} z_{in}, 2^{j_n} x) p_0(x)\, dx \\
\geq D_2^{-1} \int K(2^{j_n} z_{in}, 2^{j_n} x) \|p_0\|_\infty\, dx \\
- D_2^{-1} \int |K(2^{j_n} z_{in}, 2^{j_n} x)| \|p_0\|_\infty - p_0(x)|\, dx \\
\geq 2^{-j_n} D_2^{-1} \|p_0\|_\infty (1 - \varepsilon \|\Phi\|_1),
\end{aligned}
$$

which gives $\mu_1 = D_2^{-1} \|p_0\|_\infty (1 - \varepsilon \|\Phi\|_1)$ in (d). Third, for the lower bound in (e), note that the inequalities (27) give $(E(g_i^{(n)}(X)))^2 = O(2^{-2j_n})$, whereas

$$E(g_i^{(n)}(X))^2 = 2^{-j_n} \int_{B_1 - B_2}^{B_2 - B_1} (g_i^{(n)}(2^{j_n} z_{in}, 2^{j_n} z_{in} + u))^2 p_0(z_{in} + u 2^{-j_n})\, du$$

$$\geq 2^{-j_n} \|p_0\|_\infty (1 - \varepsilon),$$

since $z_{in} + u2^{-j_n} \in I$ and by construction of $I$. So the lower bound in condition (e) is satisfied with $\sigma_1^2 = \|p_0\|_\infty(1-2\varepsilon)$ for all $n$ large enough, which, together with $\sigma_2^2 = \|p_0\|_\infty$, gives condition (g). $\square$

This proposition, Theorem 1 and the bounds (20) already determine the a.s. rate of convergence of $\|p_n - Ep_n\|_\infty$,

$$\sqrt{\frac{n}{j_n 2^{j_n}}} \|p_n - Ep_n\|_\infty = O_{\text{a.s.}}(1) \text{ and not } o_{\text{a.s.}}(1)$$

and the same is true for the normalized quantity in Theorem 1. To obtain the exact limit (with normalization), we need the following proposition.

In the next proposition, the (at first sight somewhat awkward) condition $B_1, B_2 \in \mathbb{Z}$ is designed to include both the Haar wavelet and any continuous father wavelet with bounded support and bounded $p$-variation.

PROPOSITION 3. *Let $\phi$ be a father wavelet that satisfies Condition 2 and is uniformly continuous on $(B_1, B_2]$, where $B_1, B_2 \in \mathbb{Z}$. Suppose $P$ has a bounded uniformly continuous density $p_0$. Let $D$ be a bounded subset of $\mathbb{R}$. Then, if $j_n$ satisfies (21), we have*

$$\limsup_{n \to \infty} \sqrt{\frac{n}{(2\log 2) j_n 2^{j_n}}} \sup_{y \in D} \left| \frac{p_n(y) - Ep_n(y)}{\sqrt{\sum_k \phi^2(2^{j_n} y - k)}} \right| \leq \|p_0\|_\infty^{1/2} \quad a.s.$$

PROOF. We choose $\lambda \in (1,2)$ and $n'_k = [\lambda^k]$ (where $[a]$ denotes the integer part of $a$). Since $[\lambda^k] \leq 2[\lambda^{k-1}]$ (as $[\lambda^k]/[\lambda^{k-1}] \to \lambda < 2$) for $k$ large enough, it follows that for such $k$, the cardinality of the set $\{2^{-j_n} : n'_{k-1} < n \leq n'_k\}$ does not exceed 2 if $\tau$ in (21) equals 1, which we assume, because the proof for larger $\tau$ requires only formal changes to the present proof. Define $n_{k-1} = n'_{k-1}$ if this cardinality is 1, and otherwise let $n_{k-1}$ be the largest integer $n$ such that $j_n = j_{n'_{k-1}}$. Then we have

(28)  $$[\lambda^{k-1}] = n'_{k-1} \leq n_{k-1} < n'_k \leq n_k < n'_{k+1} = [\lambda^{k+1}]$$

and

(29)  $$j_n = j_{n'_k} = j_{n_k} \quad \text{for } n_{k-1} < n \leq n_k.$$

Let $\delta_m = 1/m$ for $m \in \mathbb{N}$. For each given $k$ and $\delta_m$, we consider the following partition of $D$. $D$ is contained in the union of at most $2 + \text{diam}(D)/(2^{-j_{n_k}}(B_2 - B_1))$ disjoint sets $(2^{-j_{n_k}}(B_1 + l), 2^{-j_{n_k}}(B_2 + l)]$, $l \in \mathbb{Z}$. Then divide each of these intervals into $m(B_2 - B_1)$ disjoint left-open right-closed subintervals $I_{k,i}$ of length $\delta_m 2^{-j_{n_k}}$ and let $z_{ki}$ be the right endpoints of the interval $I_{k,i}$ [i.e., $z_{ki} = (B_1 + m'\delta_m + l)2^{-j_{n_k}}$ for some $1 \leq m' \leq m$ and some $l \in \mathbb{Z}$]. Then the number $l_k$ of intervals $I_{ki}$ covering $D$ satisfies

$$l_k \leq 2 + \frac{\operatorname{diam}(D)}{\delta_m 2^{-j_{n_k}}} \leq \frac{c}{\delta_m 2^{-j_{n_k}}} \tag{30}$$

for some $c$ finite (and $k$ large enough depending on $m$). These intervals $I_{ki}$ also have the following property:

(31) If $z \in I_{ki}$ and $l \in \mathbb{Z}$,
then $2^{j_{n_k}} z_{ki} - l \in (B_1, B_2] \Leftrightarrow 2^{j_{n_k}} z - l \in (B_1, B_2]$,

and, for each $z_{ki}$, this happens for $B_2 - B_1$ integers $l$. As in (24) we have

$$E\bar{K}^2(2^{-j_{n_k}} z_{ki}, 2^{-j_{n_k}} X) \leq 2^{-j_{n_k}} \|p_0\|_\infty,$$

hence the maximal version of Bernstein's inequality [see Einmahl and Mason (1996)] gives that, for all $\eta > 0$,

$$\Pr\left\{ \max_{1 \leq i \leq l_k} \max_{n_{k-1} < n \leq n_k} \left| \sum_{r=1}^n (\bar{K}(2^{-j_{n_k}} z_{ki}, 2^{-j_{n_k}} X_r) - E\bar{K}(2^{-j_{n_k}} z_{ki}, 2^{-j_{n_k}} X)) \right| \right.$$

$$\left. > \sqrt{2(1+\eta)n_k 2^{-j_{n_k}} \|p_0\|_\infty \log 2^{j_{n_k}}} \right\}$$

$$\leq 2l_k \exp\left\{ -(2(1+\eta)n_k 2^{-j_{n_k}} \|p_0\|_\infty \log 2^{j_{n_k}}) \right.$$

$$\times \left( 2n_k 2^{-j_{n_k}} \|p_0\|_\infty \right.$$

$$\left. \left. + \frac{4}{3} D_1^{-1} \|\Phi\|_\infty \sqrt{2(1+\eta)n_k 2^{-j_{n_k}} \|p_0\|_\infty \log 2^{j_{n_k}}} \right)^{-1} \right\},$$

which, by the first condition in (21), is dominated by

$$2l_k \exp\left\{ -\frac{(1+\eta)\log 2^{j_{n_k}}}{1+\eta_k} \right\} \leq cm(2^{-j_{n_k}})^{(1+\eta)/(1+\eta_k)-1}$$

for some $\eta_k \to 0$ and $c$ as in (30). This is the general term of a convergent series since $\frac{1+\eta}{1+\eta_k} - 1 > 0$ and by the second condition in (21). Then Borel–Cantelli shows that

$$\limsup_n \max_{1 \leq i \leq l_k} \max_{n_{k-1} < n \leq n_k} \left( \left| \sum_{r=1}^n (\bar{K}(2^{-j_n} z_{ki}, 2^{-j_n} X_r) - E\bar{K}(2^{-j_n} z_{ki}, 2^{-j_n} X)) \right| \right)$$
(32)
$$\times (\sqrt{2n2^{-j_n} \|p_0\|_\infty \log 2^{j_n}})^{-1} \leq \lambda$$

almost surely, where we use (28) and (29).

Now consider the class of functions

$$\mathcal{G}_{ki} = \{\bar{K}(2^{j_{n_k}}z_{ki}, 2^{j_{n_k}}(\cdot)) - \bar{K}(2^{j_{n_k}}z, 2^{j_{n_k}}(\cdot)) : z \in I_{ki} \cap D\}$$

for $k \in \mathbb{N}$, $1 \leq i \leq l_k$. These classes are of VC-type—with the same VC-characteristics for each $k$ and $i$—by Lemma 2 and the permanence properties of VC-classes. We apply Talagrand's inequality to them, and we must estimate the variance $\sigma$ of the functions in these classes: by (29), (31) and change of variables, for $n_{k-1} < n \leq n_k$ and $z \in I_{ki}$ we have

$$\begin{aligned}(33)\quad & E(K(2^{j_{n_k}}z_{ki}, 2^{j_{n_k}}X) - K(2^{j_n}z, 2^{j_n}X))^2 \\ & \leq \|p_0\|_\infty (B_2 - B_1)^2 \left(\int \phi^2(x)\,dx\right) \omega^2_{\delta_m}(\phi) 2^{-j_{n_k}},\end{aligned}$$

where $\omega_\delta(\phi)$ denotes the $\delta$-modulus of continuity of $\phi$ on $(B_1, B_2]$. The same computation gives

$$\begin{aligned}(\|K(2^{j_{n_k}}z_{ki}, \cdot)\|_2 - \|K(2^{j_n}z, \cdot)\|_2)^2 &\leq \int (K(2^{j_{n_k}}z_{ki}, u) - K(2^{j_{n_k}}z, u))^2\,du \\ &\leq (B_2 - B_1)^2 \left(\int \phi^2(x)\,dx\right) \omega^2_{\delta_m}(\phi).\end{aligned}$$

Using the last two estimates, (20), (24) and that $|a/\alpha - b/\beta| \leq \alpha^{-1}|a - b| + |b|(\alpha\beta)^{-1}|\beta - \alpha|$ we obtain

$$E(\bar{K}(2^{j_{n_k}}z_{ki}, 2^{j_{n_k}}X) - \bar{K}(2^{j_n}z, 2^{j_n}X))^2 \leq C_\phi \|p_0\|_\infty (B_2 - B_1)^2 \omega^2_{\delta_m}(\phi) 2^{-j_{n_k}}$$

for $n_{k-1} < n \leq n_k$ and $z \in I_{ki}$. We set $\sigma_k^2 = C_\phi \|p_0\|_\infty (B_2 - B_1)^2 \omega^2_{\delta_m}(\phi) 2^{-j_{n_k}} := C^2 \omega^2_{\delta_m}(\phi) 2^{-j_{n_k}}$, and $U = U_k = 4D_2^{-1}\|\Phi\|_\infty$. By the first condition in (21) we have

$$\frac{n_k \sigma_k^2}{\log(U/\sigma_k)} \to \infty$$

and Talagrand's inequality (60) gives

$$\Pr\left\{\max_{1 \leq i \leq l_k} \max_{n_{k-1} < n \leq n_k} \|n(P_n - P)\|_{\mathcal{G}_{ki}} > \sqrt{3C_3 n_k \sigma_k^2 \log\frac{U}{\sigma_k}}\right\} \leq Rl_k \exp\left\{-3\log\frac{U}{\sigma_k}\right\}.$$

(Note that the supremum over $\mathcal{G}_{ki}$ is a countable supremum by Remark 2.) Now, for a fixed constant $L'$,

$$l_k \exp\left\{-3\log\frac{U}{\sigma_k}\right\} \leq cm2^{j_{n_k}}\left(\frac{\sigma_k}{U}\right)^3 \leq L'm\omega^3_{\delta_m}(\phi) 2^{-j_{n_k}/2},$$

which, by the second limit in (21) and by (28), is the general term of a convergent series. Then, for $n_{k-1} < n \leq n_k$, and $k, m$ large enough, one has

$$n_k \sigma_k^2 \log\frac{U}{\sigma_k} \leq 2\lambda^2 C^2 n 2^{-j_n} \log(2^{j_n}) \omega^2_{\delta_m}(\phi) \log\frac{U}{C\omega_{\delta_m}(\phi)}.$$

It then follows by Borel–Cantelli that

$$
\begin{aligned}
(34) \quad &\limsup_n \max_{1 \le i \le l_k} \max_{n_{k-1} < n \le n_k} \frac{\|n(P_n - P)\|_{\mathcal{G}_{ki}}}{\sqrt{n 2^{-j_n} \log 2^{j_n}}} \\
&\le \sqrt{6C_3} \lambda C \omega_{\delta_m}(\phi) \log \frac{U}{C \omega_{\delta_m}(\phi)}
\end{aligned}
$$

almost surely, for all $\lambda$ and for all $m$ large enough. Now combining (32) and (34) we have for these $m$

$$
\limsup_n \frac{\|\sum_{r=1}^n (\bar{K}(2^{-j_n}(\cdot), 2^{-j_n} X_r) - E\bar{K}(2^{-j_n}(\cdot), 2^{-j_n} X))\|_\infty}{\sqrt{2n 2^{-j_n} \|p_0\|_\infty \log 2^{j_n}}}
$$

$$
\le \lambda + \lambda C' \omega_{\delta_m}(\phi) \log \frac{U}{C \omega_{\delta_m}(\phi)},
$$

and the result follows by letting $\lambda \to 1$ and $m \to \infty$. $\square$

Summarizing, we have basically proved the main theorem of this section:

THEOREM 2. *If $\phi$, $P$ and $\{j_n\}$ are as in Proposition 3, then*

$$
\lim_{n \to \infty} \sqrt{\frac{n}{2(\log 2) j_n 2^{j_n}}} \sup_{y \in \mathbb{R}} \left| \frac{p_n(y) - E p_n(y)}{\sqrt{\sum_k \phi^2(2^{j_n} y - k)}} \right| = \|p_0\|_\infty^{1/2} \qquad a.s.
$$

PROOF. Define $D_\varepsilon = \{x \in \mathbb{R} : p_0(x) > \varepsilon, |x| < 1/\varepsilon\}$ for $\varepsilon > 0$. Applying Proposition 3 to $D_\varepsilon$ and Corollary 1 to $D_\varepsilon^c$, we obtain

$$
\limsup_{n \to \infty} \sqrt{\frac{n}{2(\log 2) j_n 2^{j_n}}} \sup_{y \in \mathbb{R}} \left| \frac{p_n(y) - E p_n(y)}{\sqrt{\sum_k \phi^2(2^{j_n} y - k)}} \right| \le \|p_0\|_\infty^{1/2} + M 2^{\tau/2} \|p_0\|_{D_\varepsilon^c}
$$

for all $\varepsilon > 0$. Now, since, $\|p_0\|_{D_\varepsilon^c} \to 0$ as $\varepsilon \to 0$ by uniform continuity of $p_0$, this lim sup does not exceed $\|p_0\|_\infty^{1/2}$. This and Proposition 2 prove the theorem. $\square$

With the natural changes in the variance computations (24) and (33), the proof of Theorem 2 implies a result similar to the one in Massiani (2003), which is the counterpart for the wavelet density estimator of the classical result of Stute (1982) for the Parzen–Rosenblatt estimator.

COROLLARY 2. *Let $\phi$ and the sequence $j_n$ be as in Proposition 3. Let $D$ be a bounded subset of $\mathbb{R}$ and assume that $p_0$ is uniformly continuous on a* *neighborhood of $D$ and $\|p_0\|_D \neq 0$. Then*

$$\lim_{n\to\infty} \sqrt{\frac{n}{2(\log 2)j_n 2^{j_n}}} \sup_{y\in D} \left|\frac{p_n(y) - Ep_n(y)}{\sqrt{\sum_k \phi^2(2^{j_n}y - k)}}\right| = \|p_0\|_D^{1/2} \quad a.s.$$

*If furthermore $\inf_{x\in D} p_0(x) > 0$, then*

$$\lim_{n\to\infty} \sqrt{\frac{n}{2(\log 2)j_n 2^{j_n}}} \sup_{y\in D} \left|\frac{p_n(y) - Ep_n(y)}{\sqrt{p_0(y)\sum_k \phi^2(2^{j_n}y - k)}}\right| = 1 \quad a.s.$$

REMARK 4 (Moments and Laplace transforms). We note that the a.s. limits in the previous theorem can be complemented by convergence of moments. In fact, direct application of inequality (61) gives that under the conditions of Theorem 1,

$$(35) \quad \sup_n E\exp\left\{\lambda\sqrt{\frac{n}{2(\log 2)j_n 2^{j_n}}} \sup_{y\in\mathbb{R}} \left|\frac{p_n(y) - Ep_n(y)}{\sqrt{\sum_k \phi^2(2^{j_n}y - k)}}\right|\right\} < \infty$$

for all $\lambda \geq 0$, and the same is true without the normalization by $\sqrt{\sum_k \phi^2(2^{j_n}y - k)}$. This yields enough uniform integrability to obtain that under the conditions of Theorem 2,

$$(36) \quad \lim_{n\to\infty} E\exp\left\{\lambda\sqrt{\frac{n}{2(\log 2)j_n 2^{j_n}}} \sup_{y\in\mathbb{R}} \left|\frac{p_n(y) - Ep_n(y)}{\sqrt{\sum_k \phi^2(2^{j_n}y - k)}}\right|\right\} = e^{\lambda\|p_0\|_\infty^{1/2}}$$

for all $\lambda \geq 0$. In particular we obtain convergence of all moments in Theorem 2 (as well as uniform boundedness of all moments in Theorem 1 and in Remark 3).

REMARK 5 (Haar wavelet and normalization). Theorem 2 (and Corollary 2) applies to the Haar father wavelet $\phi = 1_{(0,1]}$ (which is obviously uniformly continuous on $(B_1, B_2] = (0,1]$). In this case, $\sum_k \phi^2(2^{j_n}y - k) = 1$ for all $j, y$. However, if $\phi$ is not constant on $(B_1, B_2]$, the quantity $\sum_k \phi^2(2^{j_n}y - k) = \int K^2(2^{j_n}y, u)\,du$—although bounded from above and below—depends on $y$ and $j_n$, which is why it must be part of the normalization instead of the limit.

REMARK 6 (Comparison to convolution kernels). The resolution level $j_n$ in wavelet density estimation, more exactly, the quantity $2^{-j_n}$, corresponds to the window width $h_n$ in the classical ("Parzen–Rosenblatt") convolution kernel density estimator. If $K(y,x) = K(y-x)$ then the variance of the estimator $\tilde{p}_n(y) = h_n^{-1}n^{-1}\sum_{i=1}^n K((y - X_i)/h_n)$ is asymptotically of the order $n^{-1}h_n^{-1}p_0(y)\|K\|_2^2$, whereas the order of the variance of $p_n(y)$ is $n^{-1}2^{j_n}p_0(y)\int K^2(2^{j_n}y,x)\,dx$, which is different (except for Haar wavelets) since the $\mathcal{L}^2(\mathbb{R})$ norm of $K(y,x)$ oscillates with $y$. Modulo these differences, the a.s. asymptotic behavior of wavelet estimators is similar to that of convolution kernel estimators [cf. Stute (1982), Deheuvels (2000) and Giné and Guillou (2002)]. Regarding proofs, generally, the derivation of Theorems 1 and 2 follows the same pattern of proof of Theorems 2.3 and 3.3 in Giné and Guillou (2002); the short proofs of their Theorem 2.3 and of our Theorem 1 consist of a direct application of Talagrand's inequality, moment bounds for VC-type classes of functions and "blocking," and here the differences only reside in the fact that the classes of functions associated with the kernels $K$ are not the same (but in both cases of VC type), in different variance computations, and in measurability considerations. However, the proof of the exact limit law (Theorem 2) is more delicate in the wavelet case. In Proposition 3 above we cannot use continuity of translations and dilations in $\mathcal{L}^1(\mathbb{R})$ as in the upper bound part of Proposition 3.1 in Giné and Guillou (2002). Similarly, the conditions (a)–(g) that we verify in the proof of Proposition 2 also require different methods than those in Einmahl and Mason (2000) or Giné and Guillou (2002).

REMARK 7 [Nonorthogonal $\phi(\cdot-k)'s$]. The estimator $p_n(y) = P_n(K(2^{j_n}y,\cdot))$ for $K(y,x) = \sum_k \phi(y-k)\phi(x-k)$ makes sense even if $\phi$ is *not* a father wavelet, that is, the $\phi(\cdot-k)$ need not form an orthogonal system. Assuming $\phi$ satisfies Condition 2 and $\inf \|K(y,\cdot)\|_2 > 0$, then the results proved so far in this section still hold true both for $\|p_n - Ep_n\|_\infty$ and for $\sup_y |p_n(y) - Ep_n(y)|/\|K(2^{j_n}y,\cdot)\|_2$.

4.1. *Approximation error and optimal rates of convergence.* The previous results were formulated for the deviation of $p_n$ from $Ep_n$, whereas the quantity of statistical interest is $p_n - p_0$. The "bias" $Ep_n - p_0 = K_j(p_0) - p_0$ is nonrandom and can be dealt with by using standard results on approximation of functions by wavelets. If $p_0$ is uniformly continuous then, by (6), (7) and Minkowski's inequality for integrals, $\|K_j(p_0) - p_0\|_\infty \le \int \Phi(u) \|p_0(2^{-l}u + \cdot) - p_0\|_\infty \, du \to 0$ for $l \to \infty$ and $\phi$ satisfying Condition 1(0), so that if also Condition 2 holds, then

$$\|p_n - p_0\|_\infty = o_{\text{a.s.}}(1)$$

by Remark 3 if one chooses $2^{j_n} \simeq n/(\log n)^{1+\delta}$ for some $\delta > 0$.

If more is known on the smoothness of $p_0$ one can obtain rates of convergence. The approximation error in sup-norm loss of a function $f$ by wavelets is related to containment of $f$ in the Besov space $B^t_{\infty\infty}(\mathbb{R})$. Recall from (12) that these spaces include the more classical Hölder spaces $C^t(\mathbb{R})$. A function $p_0$ in $B^t_{\infty\infty}(\mathbb{R})$ is approximated by its projection $K_j(p_0)$ in the uniform norm at rate $2^{-jt}$ if $\phi$ has some regularity. To be precise, if $\phi, \psi$ satisfy Condition 1(T) for $0 < t < T + 1$, and if $p_0 \in B^t_{\infty\infty}(\mathbb{R})$, then the bounds

$$(37) \qquad \|K_j(p_0) - p_0\|_\infty \le C 2^{-jt} \quad \text{and} \quad \sup_k |\beta_{lk}(p_0)| \le C 2^{-l(t+1/2)}$$

can be shown to hold [e.g., Härdle et al. (1998), Theorem 9.4]. Inspection of the proof of that theorem shows that the constant $C$ depends on $p_0$ only through its Besov norm $\|p_0\|_{t,\infty,\infty}$.

As a consequence we have the following theorem.

THEOREM 3. *Assume that the density $p_0$ of $P$ satisfies $p_0 \in B^t_{\infty\infty}(\mathbb{R})$ for some $t > 0$. Let $p_n$ be as in (11) where $\phi$ satisfies Condition 2, and $\phi, \psi$ are such that Condition 1(T) holds for some $0 \le T < \infty$. If $j_n$ satisfies (21), then*

$$\sup_{x \in \mathbb{R}} |p_n(x) - p_0(x)| = O_{\text{a.s.}}\left(\sqrt{\frac{j_n 2^{j_n}}{n}} + 2^{-tj_n}\right)$$

*and*

$$\left(E \sup_{x \in \mathbb{R}} |p_n(x) - p_0(x)|^p\right)^{1/p} = O\left(\sqrt{\frac{j_n 2^{j_n}}{n}} + 2^{-tj_n}\right)$$

*for every $0 < t < T + 1$ and for every $1 \le p < \infty$.*

PROOF. This follows from Remarks 3 and 4 and (37). □

We note that Conditions 1(T) and 2 are satisfied for a large variety of wavelets, for example, the Haar wavelet ($T = 0$), or for sufficiently regular Daubechies wavelets (arbitrary $T \ge 0$) [cf. Härdle et al. (1998), Remark 7.1].

REMARK 8 (Optimal rates of convergence over general Besov bodies). The last theorem implies that the linear wavelet estimator with $2^{j_n} \simeq (n/\log n)^{1/(2t+1)}$ achieves the optimal [over $B^t_{\infty\infty}(\mathbb{R})$-balls] rate of convergence $((\log n)/n)^{t/(2t+1)}$ in the uniform norm for estimating $p_0$ [see, e.g., Juditsky and Lambert-Lacroix (2004) for optimality of these rates]. One might ask whether the linear wavelet estimator $p_n$ is also best possible if $p_0$ is contained in some space other than $B^t_{\infty\infty}(\mathbb{R})$, for example, in a Besov space $B^t_{pq}(\mathbb{R})$ with $t > 1/p$. The continuous (Sobolev-type) imbedding of $B^t_{pq}(\mathbb{R})$ into $B^{t-1/p}_{\infty\infty}(\mathbb{R})$ (see Remark 1) and the choice $2^{j_n} \simeq (n/\log n)^{1/(2(t-1/p)+1)}$ then give $(E\|p_n - p_0\|^r_\infty)^{1/r} = O((\log n/n)^{-(t-1/p)/(2(t-1/p)+1)})$, for all $r$, which is still the optimal rate of convergence [see, e.g., Donoho et al. (1996), Theorem 1].

**5. Uniform central limit theorems for wavelet density estimators.** Consider again the wavelet estimator $p_n(y)$ defined in (11). In this section we study the stochastic process

$$f \mapsto \sqrt{n} \int (p_n(y) - p_0(y)) f(y)\, dy = \sqrt{n}(P_n^W - P)f,$$

where $f$ varies over some Donsker class $\mathcal{F}$ of functions and where $dP_n^W(y) = p_n(y)\, dy$. Note that $P_n^W(f) = P_n(K_j(f))$. The classical case is $\mathcal{F} = \{1_{(-\infty,s]} : s \in \mathbb{R}\}$, in which case one has

$$\sqrt{n}(P_n^W - P)f = \sqrt{n}(F_n^W(s) - F(s)),$$

where $F_n^W$ and $F$ are the distribution functions of $p_n$ and $p_0$, respectively. We will show that $\sqrt{n}(P_n^W - P)$ converges in distribution in $\ell^\infty(\mathcal{F})$ to $G_P$, for various Donsker classes $\mathcal{F}$. Our proofs will in fact show $\|P_n^W - P_n\|_\mathcal{F} = o_P(1/\sqrt{n})$. The limit theorem for $P_n^W$ can then be inferred from a limit theorem for $P_n$ (using the fact that $\mathcal{F}$ is a Donsker class). In the classical case of $(F_n^W - F)$, we will also obtain a Dvoretzky–Kiefer–Wolfowitz type inequality, the compact law of the iterated logarithm, as well as a strong invariance principle.

We set, for ease of notation, $P(K_j)(y) = PK_j(y,\cdot)$, and we will use the symbol $P(K_j)$ both for the function and the finite signed measure that has it as density. The same applies to $P_n(K_j)$. For $f$, a bounded function, the following decomposition will be useful:

(38) $$(P_n^W - P_n)f = (P_n - P)(K_j(f) - f) + (P(K_j) - P)f.$$

The first term is stochastic, whereas the second ("expectation") term is deterministic, and we will deal with these two terms separately.

5.1. *CLT and strong invariance principles for the distribution function of the wavelet estimator.* We will first treat the classical special case $\mathcal{F} = \{1_{(-\infty,s]} : s \in \mathbb{R}\}$, which corresponds to studying the distribution function

$$F_n^W(s) = \int_{-\infty}^s p_n(y)\, dy$$

of the wavelet density estimator $p_n$. The key result will be an exponential inequality for the random quantity $\sqrt{n}\|F_n^W - F_n\|_\infty$, where $F_n(s) = \int_{-\infty}^s dP_n$ is the empirical distribution function. This inequality will follow from applying Talagrand's inequality (and Lemma 2) to the centered term in the decomposition (38), but we first must show that the second ("expectation") term is sufficiently small for relevant choices of $j$.

LEMMA 3. *Let $K(y,x)$ be a projection kernel as in (5) arising from the father wavelet $\phi$, and assume that $\phi, \psi$ satisfy Condition 1(T) for some* $0 \leq T < \infty$. *Assume further that the density $p_0$ is a bounded function—in which case we set $t = 0$—or that $p_0 \in B^t_{\infty\infty}(\mathbb{R})$ for some $t$, $0 < t < T + 1$. Let $\mathcal{F} = \{1_{(-\infty, s]} : s \in \mathbb{R}\}$. Then the inequality*

$$\sup_{f \in \mathcal{F}} |(P(K_j) - P)f| \leq C 2^{-j(t+1)}$$

*holds for some constant $C$ depending only on $\phi$ and $\|p_0\|_{t,\infty}$ (with $\|p_0\|_{0,\infty} = \|p_0\|_\infty$).*

PROOF. Using that the wavelet series (9) of $p_0 \in \mathcal{L}^1(\mathbb{R})$ converges in $\mathcal{L}^1(\mathbb{R})$, we have

$$K_j(p_0) - p_0 = -\sum_{l=j}^{\infty} \sum_k \beta_{lk}(p_0) \psi_{lk},$$

in the $\mathcal{L}^1$-sense. Therefore, since $f = 1_{(-\infty, s]} \in \mathcal{L}^\infty(\mathbb{R})$, we have

$$\begin{aligned}
(P(K_j) - P)f &= \int (K_j(p_0) - p_0) f \\
&= -\int \left( \sum_{l=j}^{\infty} \sum_k \beta_{lk}(p_0) \psi_{lk}(x) \right) f(x) \, dx \\
&= -\sum_{l=j}^{\infty} \sum_k \beta_{lk}(p_0) \int f(x) \psi_{lk}(x) \, dx \\
&= -\sum_{l=j}^{\infty} \sum_k \beta_{lk}(p_0) \beta_{lk}(f).
\end{aligned} \tag{39}$$

Thus, we only need to obtain bounds for the wavelet coefficients $\beta_{lk}(p_0)$ and $\beta_{lk}(f)$.

We first obtain a bound for $f$. We observe that

$$\int (K_{l+1} - K_l)(f) \psi_{lk} = \int \sum_r \beta_{lr}(f) \psi_{lr} \psi_{lk} = \beta_{lk}(f),$$

where the first identity holds by pointwise equality of the integrands, and the second because the sum has only a finite number of nonzero terms (due to compactness of the support of $\psi$) and since the $\psi_{lk}$'s are orthogonal. Therefore, we have, using also (6) with $\psi$ instead of $\phi$,

$$\begin{aligned}
\|\beta_{l(\cdot)}(f)\|_1 &\leq \int \sum_k |(K_{l+1} - K_l)(f)(x)| |\psi_{lk}(x)| \, dx \\
&\leq 2^{l/2} \left\| \sum_k |\psi(2^l(\cdot) - k)| \right\|_\infty \|K_{l+1}(f) - K_l(f)\|_1 \\
&\leq c 2^{l/2} (\|K_{l+1}(f) - f\|_1 + \|K_l(f) - f\|_1).
\end{aligned} \tag{40}$$

To bound the r.h.s., we have by Tonelli, (6) and the definition of $f$

$$\int \left| \int 2^l K(2^l y, 2^l x) f(x) \, dx - f(y) \right| dy$$

$$= \int \left| \int 2^l K(2^l y, 2^l u + 2^l y)(f(u+y) - f(y)) \, du \right| dy$$

$$\leq \int \int 2^l \Phi(2^l u) |f(u+y) - f(y)| \, du \, dy$$

$$= \int \Phi(u) \int |f(2^{-l} u + y) - f(y)| \, dy \, du$$

$$= \int \Phi(u) \left| \int_{s-2^{-l}u}^{s} dy \right| du$$

$$\leq 2^{-l} \int \Phi(u) |u| \, du.$$

Since $\Phi$ is bounded and has compact support, we conclude that

(41) $$\sup_{f \in \mathcal{F}} \|\beta_{l(\cdot)}(f)\|_1 \leq c' 2^{-l/2}$$

for some $c' \in (0, \infty)$. For the wavelet coefficients of $p_0$, we have from (37) that $\|\beta_{l(\cdot)}(p_0)\|_\infty \leq C 2^{-l(t+1/2)}$ for $0 < t < T+1$ and some constant $C$. If $t = 0$, one has the same bound since

$$|\beta_{lk}(p_0)| \leq 2^{l/2} \int |\psi(2^l x - k)| p_0(x) \, dx \leq 2^{-l/2} \|\psi\|_1 \|p_0\|_\infty$$

by a simple change of variables.

Applying these bounds to (39), we have

$$\sup_{f \in \mathcal{F}} \left| \int (K_j(p_0) - p_0) f \right| \leq \sup_{f \in \mathcal{F}} \sum_{l \geq j} \|\beta_{l(\cdot)}(p_0)\|_\infty \|\beta_{l(\cdot)}(f)\|_1$$

$$\leq c'' 2^{-j(t+1)},$$

which completes the proof. □

Using Lemmas 2 and 3 one can prove the following inequality, which is similar to Theorem 1 in Giné and Nickl (2009).

LEMMA 4. *Let $F_n(s) = \int_{-\infty}^{s} dP_n$ and $F_n^W(s) = \int_{-\infty}^{s} p_n(y) \, dy$, where $p_n$ is as in (11), $\phi$ satisfies Condition 2, and $\phi, \psi$ are such that Condition 1(T) holds for some $0 \leq T < \infty$. Assume further that the density $p_0$ of $P$ is a bounded function—in which case we set $t = 0$—or that $p_0 \in B_{\infty\infty}^{t}(\mathbb{R})$ for some $t$, $0 < t < T + 1$. Let $j$ satisfy $2^{-j} \geq d(\log n)/n$ for some $0 < d < \infty$. Then there exist finite positive constants $L := L(\|p_0\|_\infty, \phi, d)$, $\Lambda_0 :=$* $\Lambda_0(\|p_0\|_{t,\infty}, \phi, d)$ *such that for all* $n \in \mathbb{N}$ *and* $\lambda \geq \Lambda_0 \max(\sqrt{j2^{-j}}, \sqrt{n}2^{-j(t+1)})$
*we have*

$$\Pr(\sqrt{n}\|F_n^W - F_n\|_\infty > \lambda) \leq L \exp\left\{-\frac{\min(2^j\lambda^2, \sqrt{n}\lambda)}{L}\right\}.$$

PROOF. Set $\mathcal{F} = \{1_{(-\infty,s]} : s \in \mathbb{R}\}$. Using the decomposition (38) and Lemma 3 we have

$$\Pr(\sqrt{n}\|F_n^W - F_n\|_\infty > \lambda) \leq \Pr\left(n \sup_{f \in \mathcal{F}} |(P_n - P)(K_j(f) - f)| > \frac{\sqrt{n}\lambda}{2}\right)$$

by assumption on $\lambda$ (if we take $\Lambda_0 \geq 2C$), and we will apply Talagrand's inequality to the class

$$\tilde{\mathcal{F}} = \{K_j(f) - f - P(K_j(f) - f) : f \in \mathcal{F}\},$$

which is a VC-type class by Lemma 2 (and since $\mathcal{F}$ is VC)—to bound the last probability. Notice that the supremum over $f \in \mathcal{F}$ is in fact over a countable set by left continuity of the function $s \mapsto K_j(1_{(-\infty,s]}) - 1_{(-\infty,s]}$. Using the fact that $K$ is majorized by a convolution kernel $\Phi$ [cf. (6)], one establishes by the same arguments as in the proof of Theorem 1 in Giné and Nickl (2009) that $\tilde{\mathcal{F}}$ has constant envelope $U = c\|\Phi\|_1$ and that

$$\sup_{f \in \mathcal{F}} \|K_j(f) - f\|_{2,P} \leq c'2^{-j/2} =: \sigma$$

for some $0 < c' < \infty$ that depends only on $\|p_0\|_\infty$ and $\Phi$. Therefore, we have $\sigma < U/2$ and $n\sigma^2 > C \log(U/\sigma)$. Using (59) we can choose $\Lambda_0$ large enough so that

$$4^{-1}\sqrt{n}\Lambda_0\sqrt{j2^{-j}} > E$$

in the notation of Appendix, which means that $E + \sqrt{n}\lambda/4 \leq \sqrt{n}\lambda/2$. These conditions and the obvious bound $\log(1+x) \geq ((e-1)^{-1}x \vee 1)$ for $x > 0$ applied to (56) give the result. $\square$

REMARK 9 (Asymptotic equivalence of $F_n^W$ and $F_n$). Note that the variance estimate in the previous proof together with Lemma 3 [assuming $(d \log n)/n \leq 2^{-j_n}$], implies, using (57),

$$E \sup_{t \in \mathbb{R}} |F_n^W(t) - F_n(t)| \leq C\left[\sqrt{\frac{j}{2^j n}} + \sqrt{n}2^{j(t+1)}\right],$$

which is $o(1/\sqrt{n})$ if $j = j_n$ is such that $\sqrt{n}2^{j_n(t+1)} \to 0$.

The last remark suggests that—in the most interesting range of $j_n$s such as $2^{-j_n} \simeq n^{-1/(2t+1)}$—the integrated wavelet density estimator is asymptotically equivalent to the empirical distribution function (while, at the same time, delivering a reasonable estimate of the density). The exponential bound from Lemma 4 is the key to transferring several classical results for the empirical distribution function to the cdf of the wavelet density estimator, and we state below some of the more important results that can be obtained in this way.

For example, Lemma 4 implies a Dvoretzky, Kiefer and Wolfowitz (1956) type exponential bound, up to constants, for the distribution function of the wavelet estimator; namely, there exist universal constants $c_1$, $c_2$ such that for

$$\Lambda_0 \max(\sqrt{j 2^{-j}}, \sqrt{n}(2^{-j(t+1)})) \le \lambda \le \sqrt{n},$$

we have

(42) $$\Pr(\sqrt{n}\|F_n^W - F\|_\infty > \lambda) \le c_1 \exp\{-c_2 \lambda^2\}.$$

We next give the wavelet-analogue of Donsker's classical functional CLT for the empirical distribution function.

THEOREM 4. *Let $\phi, \psi$ and $p_0$ satisfy the conditions of Lemma 4 for some $t \ge 0$ and let $j_n$ satisfy $2^{-j_n} \ge d(\log n)/n$ for all $n$ and $\sqrt{n} 2^{-j_n(t+1)} \to 0$ as $n \to \infty$. If $F$ is the distribution function of $P$, then*

$$\sqrt{n}(F_n^W - F) \rightsquigarrow_{\ell^\infty(\mathbb{R})} G_P.$$

For the compact law of the iterated logarithm define

$$\mathcal{S} = \left\{ x \mapsto \int_{-\infty}^x f \, dP : \int f \, dP = 0, \int f^2 \, dP \le 1 \right\},$$

the Strassen set.

THEOREM 5. *Let $\phi, \psi$ and $p_0$ satisfy the conditions of Lemma 4 for some $t \ge 0$ and let $j_n$ satisfy $2^{-j_n} \ge d(\log n)/n$ for all $n$ and $\sup_n \sqrt{n}(2^{-j_n})^{t+1} = M < \infty$. Let $F$ be the distribution function of $P$. Then, almost surely, the sequence*

$$\left\{ \sqrt{\frac{n}{2 \log \log n}} (F_n^W - F) \right\}_{n=3}^\infty$$

*is relatively compact in $\ell^\infty(\mathbb{R})$ and its set of limit points coincides with the Strassen set $\mathcal{S}$.*

Finally, we consider the smallest admissible choice of $\lambda$ in Lemma 4. In the case $t = 0$ and the largest admissible resolution level $2^{j_n} \simeq n/\log n$, we see that we can take $\lambda \simeq \log n/\sqrt{n}$, the rate occurring in the Komlós, Major and Tusnády (1975), result on strong approximation of $\sqrt{n}(F_n - F)$ by Brownian bridges. Consequently, we have the following strong invariance principle for the integrated wavelet density estimator $F_n^W$.

THEOREM 6. *Let $\phi, \psi$ and $p_0$ satisfy the conditions of Lemma 4 with $t = 0$, and set $2^{-j_n} \simeq (\log n)/n$. Let $F$ be the distribution function of $P$. Then there exists a probability space that supports $X_1, X_2, \ldots$ i.i.d. with density $p_0$ and a sequence of Brownian bridges $B_n$ such that, for all $n \in \mathbb{N}$ and $x \in \mathbb{R}$,*

$$\Pr(\|\sqrt{n}(F_n^W - F) - B_n \circ F\|_\infty > n^{-1/2}((C + \Lambda_0')\log n + x)) \le 2n^{-2} + Me^{-\eta x},$$

*where $C, M, \eta$ are absolute constants and where $\Lambda_0' = \max(2L, \sqrt{2L}, \Lambda_0)$, with $\Lambda_0$ and $L$ as in Lemma 4. In particular, for these versions, one has*

$$\|\sqrt{n}(F_n^W - F) - B_n \circ F\|_\infty = O_{\text{a.s.}}\left(\frac{\log n}{\sqrt{n}}\right).$$

5.2. *General UCLTs for wavelet density estimators.* The question arises whether $\{1_{(-\infty, s]} : s \in \mathbb{R}\}$ in the last section can be replaced by a more general Donsker class $\mathcal{F}$. Considering the central limit theorem, such results were proved for other density estimators—such as nonparametric maximum likelihood estimators and kernel density estimators—in Nickl (2007) and Giné and Nickl (2008). We show in this section that such results can also be proved for the wavelet estimator $P_n^W$, for many classes $\mathcal{F}$, in particular for balls in general Besov spaces (hence covering Sobolev, Hölder and Lipschitz classes).

In the case of general (Besov) classes of functions, the wavelet structure will be particularly helpful, but before we turn to these classes, we show that Lemma 4 immediately implies the following result for bounded variation classes, since these are in the closed convex hull of indicator functions. A measurable function $f : \mathbb{R} \mapsto \mathbb{R}$ is of bounded variation if $v_1(f) < \infty$, cf. (13), and the class $\mathcal{F} = \{f : \|f\|_\infty + v_1(f) \le 1\}$ is a $P$-Donsker class for every $P$ [see, e.g., Dudley (1992)].

COROLLARY 3. *Let $\phi, \psi$ and $p_0$, satisfy the conditions of Lemma 4 for some $t \ge 0$. Then, if $\mathcal{F}_R = \{f \text{ right continuous} : \|f\|_\infty + v_1(f) \le 1\}$ and $L, \Lambda_0, \lambda, j$ are as in Lemma 4, we have for all $n \in \mathbb{N}$,*

$$\Pr(\sqrt{n}\|P_n^W - P_n\|_{\mathcal{F}_R} > \lambda) \le L \exp\left\{-\frac{\min(2^j \lambda^2, \sqrt{n}\lambda)}{L}\right\}.$$

If furthermore $\sqrt{n}2^{-j_n(t+1)} \to 0$ as $n \to \infty$ and if $\mathcal{F} = \{f : \|f\|_\infty + v(f)_1 \leq 1\}$, then also

$$\sqrt{n}(P_n^W - P) \rightsquigarrow_{\ell^\infty(\mathcal{F})} G_P.$$

PROOF. If $f$ is of bounded variation and right continuous, and $f(-\infty+) = 0$, then there exists a unique finite Borel measure $\mu_f$ such that $f(x) = \int 1_{(-\infty,x]}(v) \, d\mu_f(v)$. Since $(P_n^W - P_n)c = 0$ for $c$ constant [see (7)], we may assume that the elements in $\mathcal{F}_R$ all satisfy $f(-\infty+) = 0$. We then have from Fubini's theorem [using also (6)], for $f \in \mathcal{F}_R$ that $|(P_n^W - P_n)f| \leq \|F_n^W - F_n\|_\infty$. This already proves the first claim of the corollary by Lemma 4. To prove the second claim, observe that any $f \in \mathcal{F}$ is right-continuous except at most at a countable number of points, in particular there exists a right-continuous function $\tilde{f}$ such that $\tilde{f} = f$ almost everywhere. Since $P_n^W, P$ are absolutely continuous measures, we have

$$\sqrt{n}(P_n^W - P)f = \sqrt{n}(P_n^W - P)\tilde{f} = \sqrt{n}(P_n^W - P_n)\tilde{f} + \sqrt{n}(P_n - P)\tilde{f},$$

which proves the second claim by using the first and since $\mathcal{F}$ is $P$-Donsker. □

We will now prove a general central limit theorem for the wavelet density estimator, uniformly over Besov balls. The proof via the decomposition (38) necessitates that these balls be Donsker classes of functions. The following Donsker property of balls in $B_{pq}^s(\mathbb{R})$ was proved in Nickl and Pötscher (2007), and can be shown to be essentially sharp [see Nickl (2006)]. Note that under the following conditions on $s, p, q$, the Besov spaces $B_{pq}^s(\mathbb{R})$ can (and will be) viewed as spaces of bounded continuous functions.

PROPOSITION 4. *Let $\mathcal{F}$ be a bounded subset of $B_{pq}^s(\mathbb{R})$ where $1 \leq p \leq \infty$, $1 \leq q \leq \infty$, and let $P$ be a probability measure on $\mathbb{R}$. Suppose that one of the following conditions holds:*
  (a) $1 \leq p \leq 2$ and $s > 1/p$.
  (b) $2 < p \leq \infty$, $s > 1/2$, and $\int_\mathbb{R} |x|^{2\gamma} dP(x) < \infty$ for some $\gamma > 1/2 - 1/p$.
  (c) $1 \leq p < 2$, $q = 1$ and $s = 1/p$.
*Then $\mathcal{F}$ is $P$-Donsker.*

THEOREM 7. *Let $1 \leq p, q \leq \infty$ and $1/p + 1/r = 1$. Let $dP_n^W(x) = p_n(x) \, dx$ where $p_n$ is as in (11) and where $\phi, \psi$ satisfy part* (i) *of Condition 1*(T) *for some $1 \leq T < \infty$. For $0 < s < T + 1$ and $P, s, p, q$ satisfying one of the conditions in Proposition 4, let $\mathcal{F}$ be a bounded subset of $B_{pq}^s(\mathbb{R})$. Assume in addition that $p_0 \in \mathcal{L}^r(\mathbb{R})$—in which case we set $t = 0$—or that $p_0 \in B_{r\infty}^t(\mathbb{R})$ for some $t$, $0 < t < T + 1$. Suppose $\sqrt{n}2^{-j_n(t+s)} \to 0$ as $n \to \infty$. Then*

$$\sqrt{n}(P_n^W - P) \rightsquigarrow_{\ell^\infty(\mathcal{F})} G_P.$$

PROOF. We shall use throughout the proof the properties of Besov spaces summarized in Remark 1. Note first that under the conditions of the theorem, if $p = 1$, then $s > 1$ or $s = q = 1$, in which case $\mathcal{F}$ is a bounded subset of $BV(\mathbb{R})$, and the conclusion of the theorem follows from Corollary 3. So we need only consider the case $p > 1$.

We will use the decomposition (38) from above, and we first deal with the expectation term. As in (39), we obtain

$$\int (P(K_j) - P)f = \sum_{l \geq j} \sum_k \beta_{lk}(p_0)\beta_{lk}(f),$$

where one uses the conjugacy of $p$ and $r$ and the fact that the wavelet series of $p_0 \in \mathcal{L}^r(\mathbb{R})$ converges in $\mathcal{L}^r(\mathbb{R})$ if $1 \leq r < \infty$. If $t > 0$, we obtain from [Härdle et al. (1998), Theorem 9.4] that $\|\beta_{l(\cdot)}(p_0)\|_r \leq c 2^{-l(t+1/2-1/r)}$ for some finite constant $c$. In case $t = 0$ this follows from (6) and a computation similar to the one in (40), using Hölder's inequality. Similarly, it follows from the same reference, noting the obvious imbedding of $B^s_{pq}(\mathbb{R})$ into $B^s_{p\infty}(\mathbb{R})$, that we have

(43) $$\sup_{f \in \mathcal{F}} \|\beta_{l(\cdot)}(f)\|_p \leq c' 2^{-l(s+1/2-1/p)}.$$

Hence the second "expectation" term in (38) is of order

$$\sup_{f \in \mathcal{F}} \left| \int (K_{j_n}(p_0) - p_0) f \right| \leq \sup_{f \in \mathcal{F}} \sum_{l \geq j_n} \|\beta_{l(\cdot)}(p_0)\|_r \|\beta_{l(\cdot)}(f)\|_p$$

$$\leq \sum_{l \geq j_n} c'' 2^{-l(t+s+1-1/r-1/p)} \leq c''' 2^{-j_n(t+s)}$$

$$= o(1/\sqrt{n})$$

by the assumption on $j_n$.

It remains to treat the first term in (38). First observe that the class of functions

$$\bigcup_{j \geq 0} \mathcal{F}'_j := \bigcup_{j \geq 0} \{K_j(f) - f : f \in \mathcal{F}\}$$

is $P$-Donsker: by definition of the Besov norm and (10), we see that for $s'$ such that $\max(1/2, 1/p) < s' < \min(s, 1)$, $\|K_j(f)\|_{s',p,q}$ is bounded from above by $\|f\|_{s',p,q} \leq c\|f\|_{s,p,q}$, uniformly in $j$. Consequently, $\bigcup_{j \geq 0} \mathcal{F}'_j$ is contained in a ball of $B^{s'}_{pq}(\mathbb{R})$ of radius at most $c' \sup_{f \in \mathcal{F}} \|f\|_{s,p,q}$ for some constant $0 < c' < \infty$, hence it is $P$-Donsker by Proposition 4. So, in order to prove

$$\|P_n - P\|_{\mathcal{F}'_{j_n}} = o_P(1/\sqrt{n}),$$

it suffices to show that the variances $\sup_{f \in \mathcal{F}'_{j_n}} Ef^2(X)$ converge to zero. Since bounded subsets of $B^{s'}_{pq}(\mathbb{R})$ are uniformly bounded classes of functions under the conditions of the theorem, we have

$$E(K_j(f)(X) - f(X))^2 \leq c \int |K_j(f)(x) - f(x)| p_0(x)\, dx$$
(44)
$$\leq c \|K_j(f) - f\|_p \|p_0\|_r$$

and this completes the proof since $p_0 \in \mathcal{L}^r(\mathbb{R})$ by assumption and since

$$\sup_{f \in \mathcal{F}} \|K_{j_n}(f) - f\|_p \leq c' \sup_{f \in \mathcal{F}} \|K_{j_n}(f) - f\|_{s',p,q}$$

$$= \sup_{f \in \mathcal{F}} \left( \sum_{l \geq j_n}^{\infty} (2^{l(s'+1/2-1/p)} \|\beta_{l(\cdot)}(f)\|_p)^q \right)^{1/q} \to 0$$

as $n \to \infty$, by (43). □

Giné and Nickl [(2008), Theorems 5–7, Lemma 12] proved an analogue of Theorem 7 and of Corollary 3 for the classical kernel density estimator. At first sight the proof there seems somewhat more involved, but it should be noted that the proof in the wavelet case relies on nontrivial results such as the wavelet characterization of Besov spaces together with the Donsker property of Besov balls (Proposition 4), which cannot be used in the case of convolution kernels. We should also mention that the case $p \geq 2$ (and compactly supported $p_0$) in the above theorem was considered in Nickl (2007) for the much more involved case of nonparametric maximum likelihood estimators.

**6. Adaptation in sup-norm loss and the "plug-in property" of thresholding wavelet estimators.** The linear wavelet estimator $p_n(y)$ from (11) requires choosing $j_n$, and the choice of $j_n$ that produces optimal results for $p_n$ depends on the smoothness $t$ of the true density $p_0$ (cf. the discussion in Remark 8). From a practical point of view, this is a drawback, as $p_0$ is unknown. A remedy for this problem was suggested in Donoho et al. (1996) by considering so called "thresholding" wavelet estimators, defined as follows. Note first that we may write, for $j_0 < j_1$, both integers,

$$P_n(K_{j_1}) = P_n(K_{j_0}) + \sum_{l=j_0}^{j_1-1} P_n(K_{l+1} - K_l) = P_n(K_{j_0}) + \sum_{l=j_0}^{j_1-1} \sum_k \hat{\beta}_{lk} \psi_{lk}.$$

*Hard* thresholding (the only one we will consider) consists of keeping in this sum only those $\hat{\beta}_{lk}$ that are larger than a threshold $\tau$. That is, for $j_i = j_i(n)$ and $\tau = \tau(n,l)$ the hard thresholding estimator of $p_0$ is given by

$$p_n^H(y) = P_n(K_{j_0}(y,\cdot)) + \sum_{l=j_0}^{j_1-1} \sum_k \hat{\beta}_{lk} I_{\{|\hat{\beta}_{lk}|>\tau\}} \psi_{lk}(y). \tag{45}$$

It is known [e.g., Donoho et al. (1996), Juditsky and Lambert-Lacroix (2004)] that if $\tau$, $j_0$, $j_1$ are chosen in a suitable way (not requiring the knowledge of the smoothness parameter $t$), then $p_n^H$ is rate-adaptive within a logarithmic factor for any $L_{p'}$ loss, $1 \le p' < \infty$, that is

$$\sup_{p_0 \,:\, \|p_0\|_{t,p,q} \le L, |\operatorname{supp}(p_0)| \le M} E_{p_0} \|p_n^H - p_0\|_{p'}^{p'} \le C(\log n)^\gamma r_n(t, p'),$$

where $\gamma > 0$, $C$ is a constant and $r_n(t, p')$ is the minimax rate of convergence for estimating a density in the given Besov ball.

We now show, without assuming compact support for $p_0$, that the thesholding wavelet estimator is rate adaptive for supnorm loss without the logarithmic penalty and that, simultaneously, its distribution function is $\sqrt{n}$-consistent in the sup norm (in fact, it satisfies the UCLT). The pattern of proof of the result below follows that of the aforementioned authors, but we must use the results from the previous sections in several instances, and we must deal with the unbounded support of $p_0$ by introducing a moment condition for $p_0$ of arbitrarily small order combined with an application of Hoffmann–Jørgensen's inequality.

For $\kappa > 0$, define the constant

$$c(\kappa) := c(\kappa, \psi, \|p_0\|_\infty) = \frac{\kappa^2}{8\|\psi\|_2^2 \|p_0\|_\infty + 8/(3\sqrt{\log 2})\kappa \|\psi\|_\infty}.$$

Also, define

$$\mathcal{P}(L, L', \eta) = \left\{ p_0 : \|p_0\|_{t,\infty,\infty} \le L, \int |x|^\eta p_0(x)\,dx \le L' \right\}.$$

THEOREM 8. *Suppose $\phi$ satisfies Condition 2 and $\phi$, $\psi$ are such that Condition 1(T) holds for some $0 \le T < \infty$. Assume further that the density $p_0$ of $P$ satisfies $p_0 \in B_{\infty\infty}^t(\mathbb{R})$ for some $t$, $0 < t < T+1$, and that $E|X|^\eta < \infty$ for some $\eta > 0$. Let $p_n^H$, $n \ge 2$, be the thresholding estimator in (45) corresponding to*

$$\tau = \tau(n, l) = \kappa\sqrt{l/n},$$

$$2^{j_0} \simeq (n/\log n)^{1/(2(T+1)+1)} \quad \text{and} \quad n/\log n \le 2^{j_1} \le 2n/\log n, \qquad j_0 < j_1,$$

where $\kappa > 0$ is chosen so that $c(\kappa) \geq (T+3)(1+\eta^{-1})\log 2$. Then

$$\sup_{p_0 \in \mathcal{P}(L,L',\eta)} E_{p_0} \sup_{y \in \mathbb{R}} |p_n^H(y) - p_0(y)| = O\left(\left(\frac{\log n}{n}\right)^{t/(2t+1)}\right). \tag{46}$$

Moreover, letting $F_n^H$ and $F$ denote the distribution functions of $p_n^H$ and $p_0$, respectively,

$$\sqrt{n}(F_n^H - F) \rightsquigarrow_{\ell^\infty(\mathbb{R})} G_P. \tag{47}$$

PROOF. Since

$$p_0 = K_{j_0}(p_0) + \sum_{l=j_0}^{j_1-1}(K_{l+1} - K_l)(p_0) - K_{j_1}(p_0) + p_0$$

and since

$$\sum_{l=j_0}^{j_1-1}(K_{l+1} - K_l)(p_0) = \sum_{l=j_0}^{j_1-1}\sum_k \beta_{lk}(p_0)\psi_{lk}$$

with the last series converging pointwise (in fact uniformly) because $p_0 \in \mathcal{L}^1(\mathbb{R})$, we have,

$$\|p_n^H - p_0\|_\infty \leq \|(P_n - P)(K_{j_0})\|_\infty$$
$$+ \left\|\sum_{l=j_0}^{j_1-1}\sum_k (\hat{\beta}_{lk}I_{\{|\hat{\beta}_{lk}|>\tau\}} - \beta_{lk}(p_0))\psi_{lk}\right\|_\infty$$
$$+ \|K_{j_1}(p_0) - p_0\|_\infty.$$

The expectation of the first term is

$$O(((\log n)/n)^{(T+1)/(2(T+1)+1)}) = o(((\log n)/n)^{t/(2t+1)})$$

by (17) and since $t < T+1$. The third term is at most of the order $2^{-j_1 t}$ by (37), and this is $O((\log n/n)^t) = o((\log n/n)^{t/(2t+1)})$.

It remains to consider the second term, which can be decomposed as

$$\sum_{l=j_0}^{j_1-1}\sum_k (\hat{\beta}_{lk} - \beta_{lk})\psi_{lk}[I_{[|\hat{\beta}_{lk}|>\tau,|\beta_{lk}|>\tau/2]} + I_{[|\hat{\beta}_{lk}|>\tau,|\beta_{lk}|\leq\tau/2]}]$$
$$- \sum_{l=j_0}^{j_1-1}\sum_k \beta_{lk}\psi_{lk}[I_{[|\hat{\beta}_{lk}|\leq\tau,|\beta_{lk}|>2\tau]} + I_{[|\hat{\beta}_{lk}|\leq\tau,|\beta_{lk}|\leq 2\tau]}]$$
$$:= (\mathrm{I}) + (\mathrm{II}) - (\mathrm{III}) - (\mathrm{IV}),$$

where we write $\beta_{lk}$ for $\beta_{lk}(p_0)$. We first treat the "large deviations" terms (II) and (III).

For (II) we choose $\alpha \in (1, \eta+1)$ such that

$$(48) \qquad c(\kappa) \geq \frac{(T+2)\alpha}{\alpha - 1} \log 2,$$

which is possible by the condition on $c(\kappa)$, and note,

$$E \sup_x \left| \sum_{l=j_0}^{j_1-1} \sum_k (\hat{\beta}_{lk} - \beta_{lk}) \psi_{lk}(x) I_{[|\hat{\beta}_{lk}| > \tau, |\beta_{lk}| \leq \tau/2]} \right|$$

$$(49) \qquad \leq \sum_{l=j_0}^{j_1-1} 2^{l/2} \|\psi\|_\infty \sum_k [E|\hat{\beta}_{lk} - \beta_{lk}|^\alpha]^{1/\alpha}$$

$$\times [\Pr\{|\hat{\beta}_{lk}| > \tau, |\beta_{lk}| \leq \tau/2\}]^{1-1/\alpha}.$$

Now, since $\sup_x |\psi_{lk}(x)| \leq 2^{l/2} \|\psi\|_\infty$ and $E\psi_{lk}^2(X) \leq \|\psi\|_2^2 \|p_0\|_\infty = \|p_0\|_\infty$, Bernstein's inequality gives, for $l \leq j_1 - 1$ (and $n \geq e^2$),

$$\Pr\{|\hat{\beta}_{lk}| > \tau, |\beta_{lk}| \leq \tau/2\}$$

$$\leq \Pr\left\{ \left| \frac{1}{n} \sum_{i=1}^n (\psi_{lk}(X_i) - E\psi_{lk}(X)) \right| > 2^{-1} \kappa \sqrt{l/n} \right\}$$

$$(50) \qquad \leq 2 \exp\left( -\frac{\kappa^2 l}{8\|p_0\|_\infty + (8/3)\kappa \|\psi\|_\infty \sqrt{2^l l/n}} \right)$$

$$\leq 2 \exp\left( -\frac{\kappa^2 l}{8\|p_0\|_\infty + (8/(3\sqrt{\log 2}))\kappa \|\psi\|_\infty} \right)$$

$$= 2e^{-c(\kappa)l},$$

a bound which is independent of $k$. In order to estimate $\sum_k [E|\hat{\beta}_{lk} - \beta_{lk}|^\alpha]^{1/\alpha}$, we note that, by Hoffmann–Jørgensen's inequality [see the version of Corollary 1.2.7 in de la Peña and Giné (1999)], there exists a universal constant $d(\alpha)$ such that

$$(51) \qquad \|\hat{\beta}_{lk} - \beta_{lk}\|_{\alpha,P}$$

$$\leq d(\alpha) \left( \left\| \max_{1 \leq i \leq n} \left| \frac{1}{n} (\psi_{lk}(X_i) - E\psi_{lk}(X)) \right| \right\|_{\alpha,P} + \|\hat{\beta}_{lk} - \beta_{lk}\|_{1,P} \right).$$

If $\operatorname{supp} \psi \subseteq [A_1, A_2]$, we have, for the second summand,

$$\sum_k \|\hat{\beta}_{lk} - \beta_{lk}\|_{1,P} \leq 2^{l/2+1} \|\psi\|_\infty \sum_k \int_{\operatorname{supp} \psi_{lk}} dP$$

$$(52) \qquad \leq 2^{l/2+1} \|\psi\|_\infty \sum_k \int_{(A_1+k)/2^l}^{(A_2+k)/2^l} dP$$

$$\leq 2(A_2 - A_1 + 1)\|\psi\|_\infty 2^{l/2}$$

(since, for $l$ fixed, each $x \in \mathbb{R}$ is contained in at most $A_2 - A_1 + 1$ intervals $[(A_1 + k)/2^l, (A_2 + k)/2^l]$). In order to estimate the sum over $k$ of the first summands in (51), we first observe that for each $k$ and $l$, they are bounded by

$$2\left(nE\left|\frac{1}{n}\psi_{lk}(X)\right|^\alpha\right)^{1/\alpha} \leq n^{-(\alpha-1)/\alpha} 2^{(l/2)+1}\|\psi\|_\infty \left(\int_{\operatorname{supp}\psi_{lk}} dP\right)^{1/\alpha}.$$

Let $K_1 = \{k : 0 \in [(A_1 + k)/2^l, (A_2 + k)/2^l]\}$, which consists of at most $A_2 - A_1 + 1$ terms and set $K_2 = \mathbb{Z} \setminus K_1$. Then,

$$\sum_{k \in K_1} \left(\int_{\operatorname{supp}\psi_{lk}} dP\right)^{1/\alpha} \leq (A_2 - A_1 + 1)^{(\alpha+1)/\alpha} 2^{-l/\alpha}\|p_0\|_\infty^{1/\alpha}$$

and

$$\sum_{k \in K_2} \left(\int_{\operatorname{supp}\psi_{lk}} dP\right)^{1/\alpha}$$

$$\leq \sum_{k \in K_2} \frac{1}{(1 + (|A_1 + k| \wedge |A_2 + k|)/2^l)^{\eta/\alpha}} \left(\int_{(A_1+k)/2^l}^{(A_2+k)/2^l} (1 + |x|)^\eta \, dP\right)^{1/\alpha}$$

$$\leq 2^{l\eta/\alpha}\left(\sum_{k \in K_2} \frac{1}{(2^l + (|A_1 + k| \wedge |A_2 + k|))^{\eta/(\alpha-1)}}\right)^{1-1/\alpha}$$

$$\times \left(\sum_{k \in K_2} \int_{(A_1+k)/2^l}^{(A_2+k)/2^l} (1 + |x|)^\eta \, dP(x)\right)^{1/\alpha}$$

by Hölder. Since for $\lambda > 1$, $\sum_{k \in K_2} \frac{1}{(2^l + (|A_1+k| \wedge |A_2+k|))^\lambda} \leq 2\sum_{r=2^l}^\infty r^{-\lambda}$, we get

$$2^{l\eta/\alpha}\left(\sum_{k \in K_2} \frac{1}{(2^l + (|A_1 + k| \wedge |A_2 + k|))^{\eta/(\alpha-1)}}\right)^{1-1/\alpha} \leq C 2^{l(\alpha-1)/\alpha}$$

for a constant $C = C_{\eta,\alpha}$ depending only on $\eta$ and $\alpha$. Moreover,

$$\left(\sum_k \int_{(A_1+k)/2^l}^{(A_2+k)/2^l} (1 + |x|)^\eta \, dP(x)\right)^{1/\alpha}$$

$$\leq (A_2 - A_1 + 1)^{1/\alpha}(E(1 + |X|^\eta))^{1/\alpha} < \infty.$$

Thus, collecting terms,

$$\sum_k \left\|\max_{1 \leq i \leq n} \left|\frac{1}{n}(\psi_{lk}(X_i) - E\psi_{lk}(X))\right|\right\|_{\alpha,P}$$

(53)
$$\leq Cn^{-(\alpha-1)/\alpha} 2^{l/2}(2^{-l/\alpha} + 2^{l(\alpha-1)/\alpha}),$$

where $C$ depends on $\alpha$, $\eta$, $\|\psi\|_\infty$, $A_1$, $A_2$ and $\|p_0\|_\infty$. Now, adding (52) and (53) gives a bound for $\sum_k \|\hat{\beta}_{lk} - \beta_{lk}\|_{\alpha,P}$ by (51), which, combined with inequality (50), and (49), proves that the series in (49) is dominated by

$$(54) \qquad C' \sum_{l=j_0}^{j_1-1} 2^l [1 + 2(n^{-1}2^l)^{(\alpha-1)/\alpha}] e^{-c(\kappa)l(\alpha-1)/\alpha},$$

where $C'$ depends only on $\alpha$, $\eta$, $\|\psi\|_\infty$, $A_1$, $A_2$ and $\|p_0\|_\infty$. By definition of $j_1$, $n^{-1}2^l < 2/\log n$ for $l < j_1$, which gives

$$2^l[1 + 2(n^{-1}2^l)^{(\alpha-1)/\alpha}] \leq 2^l(1 + 2(2/\log n)^{(\alpha-1)/\alpha}) \leq c2^l$$

for some $c < \infty$, and, using the definition of $\alpha$ and condition (48) for $c(\kappa)$, we obtain that (54) is bounded by

$$C'' \sum_{l=j_0}^{j_1-1} 2^{-l(T+1)} \leq C''' 2^{-j_0(T+1)}$$

for suitable constants $C''$ and $C'''$. By the definition of $j_0$ and $T$, we see that this gives the bound

$$O\left(\left(\frac{\log n}{n}\right)^{(T+1)/[2(T+1)+1]}\right) = o\left(\left(\frac{\log n}{n}\right)^{t/(2t+1)}\right)$$

for the series in (49), which is what we wanted to prove for term (II).

For term (III),

$$E \sup_x \left| \sum_{l=j_0}^{j_1-1} \sum_k \beta_{lk} \psi_{lk}(x) I_{[|\hat{\beta}_{lk}| \leq \tau, |\beta_{lk}| > 2\tau]} \right|$$

$$\leq \sum_{l=j_0}^{j_1-1} 2^{l/2} \|\psi\|_\infty \sum_k |\beta_{lk}| \Pr\{|\hat{\beta}_{lk}| \leq \tau, |\beta_{lk}| > 2\tau\}$$

$$\leq c \sum_{l=j_0}^{j_1-1} 2^l e^{-c(\kappa)l} < c' 2^{-j_0(T+1)} = o\left(\left(\frac{\log n}{n}\right)^{t/(2t+1)}\right),$$

where we have used that (40) and $\|K_l(p_0)\|_1 \leq \|\Phi_l * p_0\|_1 \leq \|\Phi\|_1$ [by (6)] imply $\sum_k |\beta_{lk}| \leq C 2^{l/2}$, and that $\Pr\{|\hat{\beta}_{lk}| \leq \tau, |\beta_{lk}| > 2\tau\} \leq \Pr\{|\hat{\beta}_{lk} - \beta_{lk}| > \tau\} \leq 2\exp\{-c(\kappa)l\}$ by (50) and choice of $\kappa$.

Next, we consider (I). We will use (18) and we should note in advance that if $l \leq j_1$, then $\sqrt{l/n} \geq C 2^{l/2} l/n$, so that $\sqrt{l/n}$ is the dominating term in that bound. Let $j_1(t)$ be such that $j_0 < j_1(t) \leq j_1 - 1$ and $2^{j_1(t)} \simeq (n/\log n)^{1/(2t+1)}$ [such $j_1(t)$ exists by the definitions]. Using (18) and (6) we have

$$E \sup_x \left| \sum_{l=j_0}^{j_1(t)-1} \sum_k (\hat{\beta}_{lk} - \beta_{lk}) \psi_{lk}(x) I_{[|\hat{\beta}_{lk}| > \tau, |\beta_{lk}| > \tau/2]} \right|$$

$$\leq \sum_{l=j_0}^{j_1(t)-1} E\Big(\sup_k |\hat{\beta}_{lk} - \beta_{lk}|\Big) 2^{l/2} \sup_x \sum_k |\psi(2^l x - k)|$$

$$\leq C \sum_{l=j_0}^{j_1(t)-1} \sqrt{\frac{2^l l}{n}} = O\left(\sqrt{\frac{2^{j_1(t)} j_1(t)}{n}}\right) = O\left(\left(\frac{\log n}{n}\right)^{t/(2t+1)}\right).$$

For the second part of (I), using the same facts as in the previous computation and that $\sup_k |\beta_{lk}(p_0)| \leq c2^{-l(t+1/2)}$ by (37), we have [recall the definition of $\tau = \tau(l, n)$]

$$E \sup_x \left| \sum_{l=j_1(t)}^{j_1-1} \sum_k (\hat{\beta}_{lk} - \beta_{lk}) \psi_{lk}(x) I_{[|\hat{\beta}_{lk}|>\tau, |\beta_{lk}|>\tau/2]} \right|$$

$$\leq \sum_{l=j_1(t)}^{j_1-1} E\Big(\sup_k |\hat{\beta}_{lk} - \beta_{lk}|\Big) \frac{2}{\kappa} \sqrt{\frac{n}{l}} \sup_k |\beta_{lk}| 2^{l/2} \sup_x \sum_k |\psi(2^l x - k)|$$

$$\leq C \sum_{l=j_1(t)}^{j_1-1} 2^{-lt} = O\left(\left(\frac{\log n}{n}\right)^{t/(2t+1)}\right).$$

Finally, for term (IV), using (6) and (37) we have

$$\sup_x \left| \sum_{l=j_0}^{j_1-1} \sum_k \beta_{lk} \psi_{lk}(x) I_{[|\hat{\beta}_{lk}|\leq\tau, |\beta_{lk}|\leq 2\tau]} \right|$$

(55)
$$\leq c \sum_{l=j_0}^{j_1-1} \sup_k 2^{l/2} |\beta_{lk}| I_{[|\beta_{lk}|\leq 2\tau]}$$

$$\leq c \sum_{l=j_0}^{j_1-1} \min(2^{l/2}\delta, C2^{-lt}),$$

where $\delta = 2\kappa\sqrt{j_1/n} \geq 2\tau$ and where $C$ only depends on the Besov norm $L$ of $p_0$. To estimate this quantity, we use an idea of Donoho et al. [(1997), proof of Theorem 3, see also Delyon and Juditsky (1996)]. Set $W(l) = \min(2^{l/2}\delta, C2^{-lt})$. Clearly $\sup_{j_0 \leq l \leq j_1-1} W(l)$ is attained at $l^*$ such that $2^{l^*} = (C/\delta)^{1/(t+1/2)}$, and $W(l^*) = C^{1-r}\delta^r = C2^{-tl^*}$ where $r = t/(t+1/2)$. Hence,

$$W(l)/W(l^*) \leq \min(2^{t(l^*-l)}, 2^{l^*t+l/2}\delta/C).$$

So the last term in (55) equals

$$c \sum_{l=j_0}^{j_1-1} W(l) \leq cW(l^*) 2^{l^*t} \delta C^{-1} \sum_{j_0 \leq l < l^*} 2^{l/2} + cW(l^*) \sum_{l \geq l^*} 2^{t(l^*-l)}$$

$$\leq c'W(l^*)2^{l^*(t+1/2)}\delta C^{-1} + c'W(l^*) = c''\delta^r = O\left(\left(\frac{\log n}{n}\right)^{t/(2t+1)}\right).$$

This concludes the proof of (46).

To prove (47), observe that, with $p_n(j_1)(y) := P_n(K_{j_1}(y, \cdot))$,

$$p_n^H - p_0 = p_n(j_1) - p_0 - \sum_{l=j_0}^{j_1-1}\sum_k \hat{\beta}_{lk}\psi_{lk}I_{[|\hat{\beta}_{lk}|\leq\tau]},$$

hence the result follows from Theorem 4 since, with $\mathcal{F} = \{1_{(-\infty,s]} : s \in \mathbb{R}\}$,

$$\sup_{f\in\mathcal{F}}\left|\int\sum_{l=j_0}^{j_1-1}\sum_k \hat{\beta}_{lk}\psi_{lk}(x)I_{[|\hat{\beta}_{lk}|\leq\tau]}f(x)\,dx\right|$$

$$\leq \sup_{f\in\mathcal{F}}\sum_{l=j_0}^{j_1-1}\kappa\sqrt{l/n}\sum_k|\beta_{lk}(f)|$$

$$\leq \frac{c}{\sqrt{n}}\sum_{l=j_0}^{j_1-1}2^{-l/2}\sqrt{l} = o(1/\sqrt{n}),$$

where we use (41). □

REMARK 10 (Choice of $\kappa$). In order to choose $\kappa$ so that the constant $c(\kappa)$ satisfies the lower bound in the theorem, one needs to choose the mother wavelet $\psi$ and know a uniform bound on $\|p_0\|_\infty$. For example, if one takes the Haar basis (and hence $T = 0$), then $\|\psi\|_2 = \|\psi\|_\infty = 1$, and if one knows in addition that $\|p_0\|_\infty \leq 1$ and that a moment of $p_0$ of order one or larger exists, then the choice $\kappa = 16$ is admissible and one can adapt to the smoothness of $p_0$ up to degree $t < 1$. If no bound on $\|p_0\|_\infty$ is available, one may replace $\|p_0\|_\infty$ by $\|p_n\|_\infty$, where $p_n$ is chosen with $2^{\hat{j}_n} \simeq n/(\log n)^2$. One can then adapt arguments from Giné and Nickl (2009) to show that Theorem 8 still holds true for this (random) choice of $\kappa$.

REMARK 11 (Adaptation in the sup-norm). Aadaptive estimation of a density in sup-norm loss was considered in Tsybakov (1998) and Golubev, Lespki and Levit (2001), who worked within the framework of the Gaussian white noise model, and adapted over Sobolev balls. Considering the density model on the real line and adaptation over the (in this context) more natural classes $B_{\infty\infty}^t(\mathbb{R})$, Giné and Nickl (2009) constructed an estimator using Lepski's method that has the same properties as the hard thresholding estimator from Theorem 8 above.

## APPENDIX: TALAGRAND'S INEQUALITY AND MOMENT BOUNDS FOR VC CLASSES

Let $X_1, \ldots, X_n$ be i.i.d. with law $P$ on $\mathbb{R}$, and let $\mathcal{F}$ be a $P$-centered (i.e., $\int f \, dP = 0$ for all $f \in \mathcal{F}$) countable class of real-valued functions on $\mathbb{R}$, uniformly bounded by the constant $U$. Let $\sigma$ be any positive number such that $\sigma^2 \geq \sup_{f \in \mathcal{F}} E(f^2(X))$, set $E := E\|\sum_{i=1}^n f(X_i)\|_\mathcal{F}$ and set $V := E\|\sum_{i=1}^n f^2(X_i)\|_\mathcal{F} \leq n\sigma^2 + 16UE$ [see Talagrand (1994) for the inequality]. Then there exists a universal constant $L$ such that, for every $t \geq 0$,

$$(56) \quad \Pr\left\{\max_{k \leq n}\left\|\sum_{i=1}^k f(X_i)\right\|_\mathcal{F} \geq E + t\right\} \leq L \exp\left\{-\frac{1}{L}\frac{t}{U}\log\left(1 + \frac{tU}{V}\right)\right\}.$$

This is Talagrand's (1996) inequality, which is usually stated for $\|\sum_{i=1}^n f(X_i)\|_\mathcal{F}$ instead of for the maximum of the partial sums. However, it follows in the stated form because Talagrand's inequality can be proved [e.g., Ledoux (2001), page 144ff] by estimation of the Laplace transform of $\|\sum_{i=1}^n f(X_i)\|_\mathcal{F}$, and $\exp\{\lambda\|\sum_{i=1}^k f(X_i)\|_\mathcal{F}\}$, $k = 1, 2, \ldots$, is a submartingale, so that Doob's inequality can be applied [see also Einmahl and Mason (2000, 2005)]. We say that $\mathcal{F}$ is a VC-type class for the envelope $U$ and with VC-characteristics $A, v$ if its $\mathcal{L}^2(Q)$ covering numbers satisfy that, for all probability measures $Q$ and $\varepsilon > 0$, $N(\mathcal{F}, \mathcal{L}^2(Q), \varepsilon) \leq (AU/\varepsilon)^v$. For such classes, assuming $Pf = 0$ for $f \in \mathcal{F}$, there exists a universal constant $L'$ such that

$$(57) \quad E\left\|\sum_{i=1}^n f(X_i)\right\|_\mathcal{F} \leq L'\left(\sqrt{v}\sqrt{n\sigma^2}\sqrt{\log \frac{AU}{\sigma}} + vU \log \frac{AU}{\sigma}\right)$$

[see, e.g., Giné and Guillou (2001)]. If $\sigma < U/2$ we may replace $A$ by 1 at the price of changing the constant $L'$. Then, if

$$(58) \quad n\sigma^2 > C \log \frac{U}{\sigma}$$

for some constant $C$ we obtain

$$(59) \quad E\left\|\sum_{i=1}^n f(X_i)\right\|_\mathcal{F} \leq L''\sqrt{n\sigma^2}\sqrt{\log \frac{U}{\sigma}} \quad \text{and} \quad V \leq L'''n\sigma^2$$

for constants $L'', L'''$ that depend only on $A, v, C$. Combining these estimates with Talagrand's inequality (56), it is easy to obtain [as in Corollary 2.2 in Giné and Guillou (2002)] that there exist constants $R$ and $C_1$ depending only on $A$ and $v$ such that for all $C_2 \geq C_1$, if

$$C_1 \sqrt{n}\sigma\sqrt{\log \frac{U}{\sigma}} \leq t \leq C_2 \frac{n\sigma^2}{U}, \qquad \sigma < U/2,$$

and (58) are satisfied, then

$$\Pr\left\{\max_{k\leq n}\left\|\sum_{i=1}^{k}f(X_i)\right\|_{\mathcal{F}}\geq t\right\}\leq R\exp\left\{-\frac{1}{C_3}\frac{t^2}{n\sigma^2}\right\}, \tag{60}$$

where $C_3=\log(1+C_2/L''')/RC_2$. In particular, for $u\geq C_1$, with $\bar{L}=L'''\vee R$,

$$\Pr\left\{\max_{k\leq n}\left\|\sum_{i=1}^{k}f(X_i)\right\|_{\mathcal{F}}\geq u\sqrt{n\sigma^2}\sqrt{\log\frac{U}{\sigma}}\right\}\leq \bar{L}\exp\left\{-\frac{u\log(1+u/\bar{L})}{\bar{L}}\log\frac{U}{\sigma}\right\}.$$

These tail probabilities are of Poisson-type, and an easy (but somewhat cumbersome) computation yields that, for all $\lambda\geq 0$,

$$\begin{aligned}
E\exp&\left\{\lambda\max_{k\leq n}\frac{\|\sum_{i=1}^{k}f(X_i)\|_{\mathcal{F}}}{\sqrt{n\sigma^2}\sqrt{\log(U/\sigma)}}\right\} \\
&\leq D(A,v,C_1,\bar{L})(1+\sqrt{\lambda\bar{L}(e^{2\lambda\bar{L}/\log(U/\sigma)}-1)} \\
&\qquad\times\exp\{\lambda\bar{L}(e^{2\lambda\bar{L}/\log(U/\sigma)}-1)\}).
\end{aligned} \tag{61}$$

**Acknowledgment.** We thank two anonymous referees for a careful reading of this article and for constructive criticism that resulted in an improved exposition.

Department of Mathematics
University of Connecticut
Storrs, Connecticut 06269-3009
USA
E-mail: gine@math.uconn.edu

Statistical Laboratory
Department of Pure Mathematics
 and Mathematical Statistics
University of Cambridge
Wilberforce Road
CB3 0WB Cambridge
United Kingdom
E-mail: nickl@statslab.cam.ac.uk